\documentclass[10pt,a4paper,final]{article}
\usepackage[left=2.250cm, right=2.250cm, top=3.750cm, bottom=3.75cm]{geometry}

\usepackage{titlesec}

\titleformat*{\section}{\normalsize\bfseries}
\titleformat*{\subsection}{\normalsize\itshape} 
\titleformat*{\subsubsection}{\normalsize\itshape}  
\titleformat*{\paragraph}{\normalsize\itshape}
\titleformat*{\subparagraph}{\normalsize\itshape}
\usepackage[font={footnotesize},labelfont=bf]{caption}

\usepackage{graphics}
\usepackage{graphicx}
 \usepackage{epsfig}

\usepackage{amssymb}
 \usepackage{amsthm}
  \usepackage{amsmath}
\usepackage{tikz-cd}

 \usepackage{float}
 \usepackage{newfloat}
\usepackage{lineno}

\usepackage{ccicons}
\usepackage{epstopdf}

\usepackage{multirow}
\usepackage{color}
\usepackage{hhline}
\usepackage{booktabs}

\definecolor{nverde}{RGB}{0,61,0} 
\definecolor{cr1}{RGB}{200,0,0}
\definecolor{cr2}{RGB}{0,0,200}
\definecolor{cr12}{RGB}{100,0,100}

\graphicspath{{./}{./figures/}}

\usepackage{authblk}
\usepackage{fancyhdr}
\usepackage{changepage}
\newcommand{\address}[2]{\affil[#1]{#2}}
\newcommand{\corref}{\footnote{Corresponding author}}
\newcommand{\MPI}{1}

 \makeatletter
\def\blfootnote{\gdef\@thefnmark{}\@footnotetext}
\makeatother

\title{Mass, momentum and energy preserving 
FEEC and broken-FEEC schemes for the incompressible Navier-Stokes equations} 

\author[\MPI]{\underline{Valentin Carlier}\corref} 
\author[\MPI]{Martin Campos Pinto}
\author[\MPI]{Francesco Fambri}
\address{\MPI}{Max-Planck-Institut f\"ur Plasmaphysik, Boltzmannstra{\ss}e 2, D-85748 Garching, Germany}

\date{22 June, 2023}

\newcommand{\keywords}[1]{\textbf{keywords}-- #1}

\allowdisplaybreaks

\usepackage[utf8]{inputenc}
\usepackage{geometry}

\usepackage[bbgreekl]{mathbbol}
\usepackage{amsmath,amssymb}
\usepackage{amsfonts}
\usepackage{bm}
\usepackage{graphicx}
\usepackage{caption}
\usepackage{subcaption}
\usepackage{authblk}
\usepackage[all]{xy}
\usepackage{tikz-cd}
\usepackage{hyperref}
\usepackage{cleveref}
\usepackage{empheq,etoolbox}
\usepackage{showlabels}
\usepackage{comment}

\usepackage{scalerel}
\usetikzlibrary{svg.path}
\definecolor{orcidlogocol}{HTML}{A6CE39}
\tikzset{
  orcidlogo/.pic={
    \fill[orcidlogocol] svg{M256,128c0,70.7-57.3,128-128,128C57.3,256,0,198.7,0,128C0,57.3,57.3,0,128,0C198.7,0,256,57.3,256,128z};
    \fill[white] svg{M86.3,186.2H70.9V79.1h15.4v48.4V186.2z}
                 svg{M108.9,79.1h41.6c39.6,0,57,28.3,57,53.6c0,27.5-21.5,53.6-56.8,53.6h-41.8V79.1z M124.3,172.4h24.5c34.9,0,42.9-26.5,42.9-39.7c0-21.5-13.7-39.7-43.7-39.7h-23.7V172.4z}
                 svg{M88.7,56.8c0,5.5-4.5,10.1-10.1,10.1c-5.6,0-10.1-4.6-10.1-10.1c0-5.6,4.5-10.1,10.1-10.1C84.2,46.7,88.7,51.3,88.7,56.8z};
  }
}
\newcommand\orcidicon[1]{\href{https://orcid.org/#1}{\mbox{\scalerel*{
\begin{tikzpicture}[yscale=-1,transform shape]
\pic{orcidlogo};
\end{tikzpicture}
}{|}}}}

\patchcmd{\subequations}
  {\theparentequation\alph{equation}}
  {\theparentequation.\alph{equation}}
  {}{}

\newtheorem{theorem}{Theorem}
\newtheorem{proposition}[theorem]{Proposition}%
\newtheorem{assumption}[theorem]{Assumption}%
\newtheorem{remark}{Remark}%

\newtheorem{lemma}[theorem]{Lemma}

\numberwithin{equation}{section}

\definecolor{darkgreen}{rgb}{.1,.5,0}

\providecommand{\martin}[1]{\textcolor{blue}{(Martin: #1)}}

\def\ve{{\varepsilon}}

\def\dd{\, {\rm d}} 
\def\Dt{{\Delta t}}

\providecommand{\norm}[1]{\lVert#1\rVert}

\def\tb{\hbox{$\|\kern -.09em |$}}
\def\bigtb{\hbox{$\big\|\kern -.09em \big|$}}
\def\Bigtb{\hbox{$\Big\|\kern -.09em \Big|$}}

\usepackage{amsopn}

\let\div\relax   
\DeclareMathOperator{\div}{div}
\DeclareMathOperator{\Div}{div}

\DeclareMathOperator{\curl}{curl}

\DeclareMathOperator{\grad}{grad}

\DeclareMathOperator{\Tcurl}{\widetilde{\curl}}

\DeclareMathOperator{\Bcurl}{\overline{\text{curl}}}
\DeclareMathOperator{\Tgrad}{\widetilde \grad}
\DeclareMathOperator{\Bgrad}{\overline{\grad}}

\usepackage{mathrsfs}

\newcommand{\iTilde}{\Tilde{\imath}}

\newcommand{\ff}{{\mathbf f}}

\newcommand{\uu}{{\mathbf u}}

\newcommand{\vv}{{\mathbf v}}
\newcommand{\ww}{{\mathbf w}}

\newcommand{\nn}{{\mathbf n}}
\newcommand{\ab}{{\mathbf a}}
\newcommand{\ee}{{\mathbf e}}
\newcommand{\Po}{\mathcal{P}}
\newcommand{\R}{{\mathbf R}}

\newcommand{\buu}{\Bar{\uu}}
\newcommand{\uue}{\mathring{\uu}}

\newcommand{\hK}{\hat{K}}
\newcommand{\hV}{\hat{V}}
\newcommand{\huu}{\hat{\uu}}
\newcommand{\hvv}{\hat{\vv}}
\newcommand{\hww}{\hat{\ww}}

\begin{document} 
\maketitle 
 
\begin{abstract}
In this article we propose two finite element schemes for the Navier-Stokes equations, based on a reformulation that involves differential operators from the de Rham sequence and an advection operator with explicit skew-symmetry in weak form.
Our first scheme is obtained by discretizing this formulation with conforming
FEEC (Finite Element Exterior Calculus) spaces: it preserves the pointwise divergence free constraint of the velocity, its total momentum and its energy, in addition to 
being pressure robust.
Following the broken-FEEC approach, our second scheme uses fully discontinuous spaces and local conforming projections to define the discrete differential operators. It preserves the same invariants up to a dissipation of energy to stabilize numerical discontinuities.
For both schemes we use a middle point time discretization which preserve these invariants at the fully discrete level and we analyse its well-posedness in terms of a CFL condition. Numerical test cases performed with spline finite elements allow us to verify the high order accuracy of the resulting numerical methods, as well as their ability to handle general boundary conditions.
\end{abstract}

\keywords{incompressible Navier-Stokes; symmetry-preserving; discrete invariants; FEEC; broken-FEEC.}	

\paragraph{Orcid numbers.} 
V. Carlier \orcidicon{0000-0002-8228-171X},
M. Campos Pinto \orcidicon{0000-0002-8915-1627},
F. Fambri \orcidicon{0000-0002-6070-8372}\blfootnote{
\underline{\textit{valentin.carlier@ipp.mpg.de}} (V. Carlier),  
\textit{martin.campos-pinto@ipp.mpg.de} (M. Campos-Pinto),  
\textit{francesco.fambri@ipp.mpg.de} (F. Fambri).  
}.

\section{Introduction}\label{sec:intro}
  
The preservation of physical invariants and underlying geometrical structure is known to play a crucial role in the simulation of various physical systems. Building a scheme that preserves such quantities allows for more stable simulations, which is of primal importance when one wants to perform long time simulations, for applications such as fluid modeling \cite{arakawa1997computational}, weather forecast \cite{phillips1959example} or fusion plasma physics \cite{lewis1970energy}.

The geometrical structure of the equations is often better understood when written in the language of differential geometry, using differential forms to rewrite various objects of the Maxwell equations \cite{bossavit1988whitney,monk2003finite}, or the diffeomorphism group for fluid evolution \cite{arnold2008topological}. 
Much effort has been put into the investigation of discretizations that preserve those underlying differential structures \cite{Bossavit.1998.ap,pavlov2011structure,gawlik2011geometric}. 
The construction of discretizations that preserve vector calculus identities 
such as $\curl \grad= 0$ and $\div \curl = 0$, or the more general exterior calculus identity $d \circ d = 0$ \cite{hiptmair2002finite,hirani2003discrete} allows the preservation of some important physical invariants. Among those discretizations, 
the FEEC framework \cite{arnold_falk_winther_2006_anum,arnold2018finite},
which provides high order discretization of the de Rham complex, has been used for a variety of application such as the Poisson equation in mixed form \cite{alonso1996error}, 
the Maxwell equations \cite{bossavit1988whitney,hiptmair2002finite} 
or the Vlasov-Maxwell system \cite{kraus2017gempic}. 
This framework was recently extended to broken spaces \cite{campos_pinto_gauss-compatible_2016,campos2022broken,gucclu2022broken}, leading to discrete operators with better locality properties and potentially more efficient discretizations.

Application of these geometric discretizations to fluid equations has provided some efficient schemes with excellent conservation properties, based on a vorticity reformulation of the Navier-Stokes system \cite{elcott2007stable,mohamed2016discrete}. However, the latter do not preserve invariants such as energy, which is known to play a crucial role in long time simulations \cite{arakawa1997computational}. Energy preserving discretizations based on the framework of Arnold \cite{arnold2008topological} have been provided in a low order discrete exterior calculus framework \cite{pavlov2011structure,gawlik2011geometric} and extended to high order in \cite{natale2018variational,gawlik2021variational}. Those geometric schemes are nevertheless not based on FEEC spaces which are particularly well adapted to
magneto-hydrodynamics (MHD) extensions in plasma modelling.
Attempts to use the FEEC framework for the simulation of fluids originate from geophysical simulations, first for the simulation of shallow-water equations \cite{cotter2012mixed,cotter2014finite} and then extended to compressible Euler equations \cite{natale2016compatible}. 
Recently it was used in \cite{palha_mass_2017} with a coupled velocity/vorticity formulation, to provide a high-order (spectral) and highly conservative scheme for the 2D incompressible Navier-Stokes equations.
This approach was extended in \cite{zhang2022mass} where a dual formulation was proposed, involving two approximations for both the velocity and the vorticity,
in order to preserve both the energy and helicity at the discrete level.

Another attractive feature of FEEC schemes for incompressible flows is their natural preservation of Helmholtz decompositions at the discrete level, which is a key ingredient
for the pressure robustness \cite{linke_role_2014,John_Linke_Merdon_Neilan_Rebholz_2017}.
We may quote for instance \cite{hanot_arbitrary_2022} where a conservative 
FEEC discretization is studied, with proven convergence and pressure-robustness properties.

As previously mentionned, using broken spaces instead of fully continuous finite element methods, allows for sparser and more local operators, leading to more efficient numerical schemes where the highly parallel architecture of modern high performance computing facilities can be better exploited. One of the most important families of numerical schemes that is built on broken finite-element spaces is surely the discontinuous Galerkin (DG) method. 
While difference schemes \cite{HarlowWelch,Patankar1972,patankar,vanKan1986} and finite-elements methods \cite{Taylor1973,Brooks1982,Hughes1986,Fortin1981,Verfuerth91,Heywood1982,Heywood1988,Arnold1984,Brezzi1989}  
have long been the state of the art, in the last decades discontinuous Galerkin methods have also been adapted to solve the low-Mach or incompressible Navier-Stokes equations.
Indeed, the design of stable DG schemes for the incompressible Navier-Stokes equations is not a trivial task. 
From one hand, physical viscosity is typically introduced in the PDE by means of second- (parabolic) or high-order  spatial derivatives, and, for DG schemes, this may introduce several difficulties, see \cite{CockburnShu1998,yan2002}. 

On the other hand, in the low-Mach limit, compressible gas-dynamics becomes asymptotically elliptic, and implicit solvers are one of the most effective powerful tool to build stable schemes. 
One of the main advantages of DG schemes with respect to standard finite-elements is the \emph{locality}, without compromises in terms of order of accuracy. Indeed, even the inversion of a mass matrix may introduce a lot of difficulties in terms for computational efficiency and related memory consumption for standard finite-element schemes. 
However, the mass matrix of DG schemes is element-local and its inversion is easy, if not eventually trivial.
To deal with the elliptic coupling between velocity and pressure, the family of fractional-step  or projection methods is probably representative of the most effective strategy to solve the incompressible Navier-Stokes equations. The main idea is to decouple the original saddle-point problem, in a sequence of decoupled elliptic problems for which the computational costs are drastically reduced. For an overview of projection methods, see \cite{GUERMOND2006}.
To this category belongs any pressure- or velocity- correction methods, as well as  standard
semi-implicit methods.

In the last decade, an interesting class of efficient symmetry-preserving DG schemes for incompressible flows has been designed on structured, adaptive staggered grids \cite{DumbserCasulli2013,FambriDumbser,AMRDGSI,Fambri2020} 
as well as unstructured grids in two and three-space dimensions \cite{TavelliDumbser2014b,TavelliDumbser2015,TavelliDumbser2016}. 
Preserving the (skew) symmetry of the differential operators has a double benefit: first, symmetric system are solvable by means of efficient methods, such as matrix-free Conjugate Gradient methods; second, in the evaluation of the discrete energy balance, symmetric terms cancels each other improving the conservation properties.
 Moreover, the theory of  generalized locally Toeplitz Algebra has also been used as a tool to describe the spectral features of the generated matrices, providing a fertile ground to the design of novel effective preconditioners \cite{GLT2018DG,Mazza2022,Barbarino2022}.

On non-staggered grids, many alternative formulations have been designed combining  the efficiency and stability properties of semi-implicit methods within the high-order framework of DG schemes, e.g. see \cite{Dolejsi2004,Dolejsi2007,Dolejsi2008} for compressible convection-diffusion flows equations, or \cite{GiraldoRestelli,Tumolo2013} for applications to the shallow water equations.
We should cite also a completely different but valuable approach that aims, instead, at relaxing the strong coupling between velocity and pressure without compromising the consistency with the incompressible limit, see e.g. \cite{bassi_artificial_2018}.  In its original version, see \cite{ElsworthToro92,Bassi2006},   an artificial compressibility perturbation is added only at the level of element-interface fluxes and, in contrast to standard artificial compressibility methods, this is done without adding pressure time derivative to the divergence-free equation.  In this direction, a novel interesting idea was proposed by \cite{Massa2022} where an exact Riemann-solver based on an artificial equation of state is designed  so that the PDE remains first-order and hyperbolic, while the incompressible Euler equations are recovered in the infinite artificial sound-speed limit.

In this work we propose new FEEC and broken-FEEC discretizations of the incompressible Navier-Stokes equation
which preserve both the energy and momentum, in addition to being conservative and pressure-robust. 
Our approach is based on a skew-symmetric reformulation of the advection operator detailed in \cref{sec:reformulation_INS}. Structure preserving discretizations are presented in \cref{sec:discretizations}:
first with a conforming FEEC scheme, which is shown to preserve mass, momentum and energy. 
Then a broken-FEEC variant is presented which preserves the mass and momentum conserving properties, and dissipates energy to stabilize the non-conforming errors. 
Both schemes are shown to be pressure robust.
\Cref{sec:Time_discretization} then presents a time discretization that is shown to preserve the invariants in the fully discrete setting. It involves an iterative procedure which is analyzed in \cref{sec:CFL,sec:scale_CFL} where a CFL condition is derived. The discrete modeling of boundary conditions is described in \cref{sec:BC}. Numerical results are finally shown in \cref{sec:Numerics} where several test-cases demonstrate the high order accuracy and conservation properties of our scheme. 
\Cref{sec:Conclusion} gathers some concluding remarks.

\section{Reformulation of the Navier-Stokes equation}
\label{sec:reformulation_INS}
We are interested in the discretization of the incompressible Navier-Stokes equation:
\begin{equation}
    \label{eqn:NS}
    \frac{\partial \uu}{\partial t}+ (\uu \cdot \grad) \uu + \grad p - \nu \Delta \uu=0 ~ ,
\end{equation}
the incompressibility condition being 
\begin{equation}
    \label{eqn:incomp}
    \div \uu =0 ~ .
\end{equation}
Here $\uu$ is the velocity vector field and $p$ the pressure. 
The equation is posed in a domain $\Omega \subset \R^d$, $d \in \{2,3\}$, 
assuming first no or periodic boundaries.
Throughout this article we will focus on the two-dimensional case,
as our numerical simulations are also conducted in this case.
However our schemes readily extend to the three-dimensional case.

\subsection{Rewriting in terms of de Rham operators}

To rewrite these equations in the framework of differential geometry,
following \cite{arnold_falk_winther_2006_anum} we a consider a primal de Rham sequence
\[
\xymatrix{
 V^0 \ar[r]^-{\curl} & V^1 \ar[r]^-{\Div} & V^2
}
\]
with $V^0 \subset H(\curl,\Omega)$, $V^1 \subset H(\Div;\Omega)$ and $V^2 \subset L^2(\Omega)$,
and a dual sequence 
\[
\xymatrix{
  \Tilde{V}^0 \ar[r]^-{\Tgrad} & \Tilde{V}^1 \ar[r]^-{\Tcurl} & \Tilde{V}^2~.
}
\]
Here the primal curl is the vector valued $\curl q = (\partial_2 q, -\partial_1 q)$
and the dual differential operators $\Tgrad$ and $\Tcurl$ are the adjoint of 
the primal ones, in the sense that 
\begin{equation} \label{eqn:dual}
  \int_\Omega \vv\cdot \Tgrad q  = -\int_\Omega q \Div \vv
  \qquad \text{ and } \qquad
  \int_\Omega \omega \Tcurl \ww  = \int_\Omega \ww \cdot \curl\omega
\end{equation}
hold for all $q \in \Tilde{V}^0, \vv \in V^1$,
and all $\ww \in \Tilde{V}^1, \omega \in V^0$.
At the continuous level they coincide with the usual gradient and scalar-valued 
curl operators, with dual spaces 
$\Tilde{V}^0 \subset H^1(\Omega)$, $\Tilde{V}^1 \subset H(\curl;\Omega)$ 
and $\Tilde{V}^2 \subset L^2(\Omega)$ equipped with proper (e.g. periodic)
boundary conditions. At the discrete level however, the dual 
sequence will be approximated in a weak form.
We will also need the interior product
\begin{equation} \label{i0}
  \iTilde_{\ab}^0 : \Tilde{V}^1 \longrightarrow \Tilde{V}^0 ~ , \quad \vv \mapsto \vv \cdot \ab
\end{equation}
associated with a given vector field $\ab$.

We first consider a weak formulation
of the Navier-Stokes \cref{eqn:NS}
where $\uu \in (H^1(\Omega))^d$ and $p \in \Tilde{V}^0$ 
(with implicit dependence with respect to time)
must satisfy
\begin{equation}
    \label{eqn:weakNS}
    \int_\Omega \frac{\partial \uu}{\partial t} \cdot \vv+ \int_\Omega(\uu \cdot \grad) \uu \cdot \vv + \int_\Omega \Tgrad p \cdot \vv + \nu \int_\Omega \Tcurl \uu \ \Tcurl \vv=0
\end{equation}
for any test function $\vv \in (H^1(\Omega))^d$, see e.g. \cite{temam_navier-stokes_1977}. 
Here we have used $-\Delta \uu = \curl \Tcurl \uu + \Tgrad \Div \uu$ and \eqref{eqn:incomp}.
Using an integration by parts, together with the divergence free condition and the absence of boundary terms, we then rewrite 
\[
    \int_\Omega(\uu \cdot \grad) \uu \cdot \vv 
    := \frac{1}{2}\Big(\int_\Omega(\uu \cdot \grad) \uu \cdot \vv -\int_\Omega(\uu \cdot \grad) \vv \cdot \uu \Big)  
    = \int_\Omega \uu \cdot s(\uu,\vv) ~ ,
\]
with 
\begin{equation} \label{eqn:def_s}
  s(\uu,\vv) := \frac{1}{2} \big(\grad \uu \cdot \vv - \grad \vv \cdot \uu\big) = \frac 12 \Big(\sum\limits_{k=1}^{d} \iTilde^0_{\ee_k} \vv \Tgrad \iTilde^0_{\ee_k} \uu - \iTilde^0_{\ee_k} \uu \Tgrad \iTilde^0_{\ee_k} \vv\Big)~.
\end{equation}
Here the $\ee_k$'s are the canonical basis vectors in $\R^d$.

This skew-symmetric convection operator will play a key role in our approach, in particular
after discretization it will naturally lead to a numerical scheme 
that preserves energy and momentum. 

Using \eqref{eqn:def_s} we rewrite the system as: 
find $\uu \in (H^1(\Omega))^d$ 
and $p \in \Tilde{V}^0$ such that 

\begin{subequations}
\label{eqn:DeRhamINS}

  \begin{empheq}[left=\empheqlbrace]{align}
    \label{eqn:DeRhamINSmom}
    \int_\Omega \frac{\partial \uu}{\partial t} \cdot \vv+ \int_\Omega \uu \cdot s(\uu,\vv) + \int_\Omega \Tgrad p \cdot \vv + \nu \int_\Omega \Tcurl \uu \ \Tcurl \vv=0 &, \qquad \forall \vv \in (H^1(\Omega))^d \\ 
    \label{eqn:DeRhamINSdiv}
    \Div \uu =0 & .
  \end{empheq}
\end{subequations}

\subsection{Conservation properties}
\label{sec:consprop}

We now recall some well known preservation properties of the Incompressible Navier-Stokes system,
in addition to the divergence-free constraint \eqref{eqn:incomp} 
which corresponds to a mass conservation property.
Even though they are common we detail their proofs based on our reformulation. 
This will emphasize the key properties of the operators leading to these invariants:
these properties will then be preserved in both our conforming and non-conforming discretizations. 

\begin{proposition}[Preservation of momentum]
\label{prop:presmom}
Solutions to \cref{eqn:DeRhamINS} preserve the momentum: 
\begin{equation}
    \label{eqn:momentconv}
    \frac{\partial}{\partial t}\int_\Omega \uu =0 ~ .
\end{equation}
\end{proposition}

\begin{proof}[Proof.]
  
Test \cref{eqn:DeRhamINSmom} with $\vv=\ee_i$: this gives 
\begin{equation}
\label{eqn:evolmom}
    \frac{\partial}{\partial t}\int_\Omega \uu_i  = -\int_\Omega \uu \cdot s(\uu,\ee_i) - \int_\Omega \Tgrad p \cdot \ee_i - \nu \int_\Omega \Tcurl \uu \ \Tcurl \ee_i  ~ .
\end{equation}
Since $\iTilde^0_{\ee_k} \ee_i = \delta_{i,k}$ we have 
$s(\uu,\ee_i)= \frac 12(\Tgrad \iTilde^0_{\ee_i} \uu - \iTilde^0_{\ee_i} \uu \Tgrad 1)$.
The duality relations \eqref{eqn:dual} characterizing $\Tgrad 1$ yield then
$\int_\Omega \Tgrad 1 \cdot \vv = -\int_\Omega 1 \Div \vv =0$ for all $\vv \in V^1$, so $\Tgrad 1=0$.
Using \cref{eqn:DeRhamINSdiv} we next compute
\[
\int_\Omega \uu \cdot s(\uu,\ee_i)= \frac 12\int_\Omega \uu \cdot \Tgrad \iTilde^0_{\ee_i} \uu 
  = -\frac 12\int_\Omega (\Div \uu) \iTilde^0_{\ee_i} \uu = 0
\]
and we also  have $\int_\Omega \Tgrad p \cdot \ee_i = -\int_\Omega p \Div \ee_i = 0$.
Using finally the duality relation \eqref{eqn:dual} for the operator $\Tcurl$ 
we compute 
$\int_\Omega \Tcurl \ee_i \cdot \omega = \int_\Omega \ee_i \cdot \curl \omega 
= \int_\Omega (\partial_{i+1} \omega_{i-1} - \partial_{i-1} \omega_{i+1}) = 0$ 
since we have no boundaries (using a cyclic notation for indices: $\omega_{d+1}=\omega_1$ and
$\partial_{d+1}=\partial_1$). 
Thus, all three terms in \cref{eqn:evolmom} are zero, which proves the claim.
\hfill\end{proof}

\begin{proposition}[Preservation of energy]
\label{prop:presener}

If $\nu=0$, solutions to \cref{eqn:DeRhamINS} preserves the (kynetic) energy: 
\begin{equation}
    \label{eqn:energyconv}
    \frac{\partial}{\partial t}\int_\Omega \frac{1}{2}||\uu||^2 =0 ~ .
\end{equation}
\end{proposition}

\begin{proof}[Proof.]
Using \cref{eqn:DeRhamINSmom} we have 
\begin{equation}
\label{eqn:evolener}
    \frac{\dd}{\dd t}\int_\Omega \frac{1}{2}||\uu||^2 
      = \int_\Omega \frac{\partial \uu} {\partial t} \cdot \uu 
      = -\int_\Omega \uu \cdot s(\uu,\uu) - \int_\Omega \Tgrad p \cdot \uu - \nu \int_\Omega \Tcurl \uu \ \Tcurl \uu ~ .
\end{equation}
The first term is 0 due to the skew symmetry of $s$.
For the second one we use again \eqref{eqn:dual} which gives
$\int_\Omega \Tgrad p \cdot \uu = - \int_\Omega p \Div \uu=0$,
and the third one vanishes for $\nu=0$. 
\hfill\end{proof}

\begin{remark}
Using the same argument as in the previous proposition we see that in the viscous case 
the energy vanishes at rate $-\nu \int_\Omega \Tcurl \uu \ \Tcurl \uu$. 
Having a correct dissipation rate for $\nu \neq 0$ is also important for 
numerical schemes aiming for long time simulations.
\end{remark}

We finally remind a important invariance property, extensively discussed in 
\cite{John_Linke_Merdon_Neilan_Rebholz_2017}. 
\begin{proposition}[Pressure-robustness property]
\label{prop:fip}
For the problem with forcing term 
\begin{equation} \label{eqn:NS_f}
  \frac{\partial \uu}{\partial t}+ (\uu \cdot \grad) \uu + \grad p - \nu \Delta \uu = \ff
\end{equation} 
the following invariance property holds:
\begin{equation} \label{fip}
  \ff \to \ff + \grad \psi  \quad \implies \quad 
  (\uu, p) \to  (\uu, p+\psi).
\end{equation} 
\end{proposition}
\begin{proof}[Proof.]
The verification is straightforward.
\hfill\end{proof}

\section{Structure preserving spatial discretization}
\label{sec:discretizations}

\subsection{Discrete De Rham sequence}

To discretize \cref{eqn:NS} we first use conforming spaces $V^{1,c}_h \subset V^1$ and $V^{2,c}_h \subset V^2$ and suppose the existence of projections $\Pi_0$, $\Pi_1$, $\Pi_2$ such that the following diagram is commuting: 
\begin{equation} \label{Vc}
  \xymatrix{
   V^0 \ar[d]_-{\Pi_0} \ar[r]^-{\curl} & V^1 \ar[d]_-{\Pi_1} \ar[r]^-{\Div} & V^2 \ar[d]_-{\Pi_2}\\
   V^{0,c}_h \ar[r]_-{\curl} & V^{1,c}_h \ar[r]_-{\Div} & V^{2,c}_h 
  }
\end{equation}
Following the standard FEEC approach these finite dimensional spaces 
will be identified with their duals. 
Assuming first periodic or homogeneous boundary conditions, 
we define the dual discrete gradient and curl operators
\begin{equation} \label{eqn:dualh}
  \Tgrad_h: V^{2,c}_h \to V^{1,c}_h
  \qquad \text{ and } \qquad
  \Tcurl_h: V^{1,c}_h \to V^{0,c}_h
\end{equation}
as the adjoints of (minus) the divergence and the curl operators, respectively.
Thus, we set
\begin{equation} \label{eqn:dualh_eq}
  \int_\Omega \Tgrad_h q_h \cdot \vv_h = -\int_\Omega q_h \Div \vv_h
  \qquad \text{ and } \qquad
  \int_\Omega (\Tcurl_h \vv_h) \omega_h = \int_\Omega \vv_h \cdot \curl\omega_h
\end{equation}
for all $q_h \in V^{2,c}_h$, $\vv_h \in V^{1,c}_h$ and $\omega_h \in V^{0,c}_h$.
From the property $\Div \curl = 0$ satisfied by the strong operators, one easily verifies that  
the discrete adjoint operators satisfy
\begin{equation} \label{dual_seq_prop}
  \Tcurl_h \Tgrad_h q_h = 0 \quad \forall q_h \in V^{2,c}_h~.
\end{equation}
We further define a discrete interior product approximating \eqref{i0} as 
\begin{equation} \label{i0h}
  \iTilde^{0h}_{\ab} \uu_h := \Tilde{P}_h^0(\ab \cdot \uu_h)
\end{equation}
where $\Tilde{P}_h^0$ is the $L^2$ projection on $V^{2,c}_h$ 
(which also discretizes the space $\Tilde{V}^0$).

This allows us to define a discrete version of $s$: 
\begin{equation}
\label{eqn:def_sh}
s_h(\uu_h,\vv_h) := \frac{1}{2} P_h^1 \Big( \sum\limits_{k=1}^{d} (\iTilde^{0h}_{\ee_k} \vv_h )\Tgrad_h (\iTilde^{0h}_{\ee_k} \uu_h )
    - (\iTilde^{0h}_{\ee_k} \uu_h) \Tgrad_h (\iTilde^{0h}_{\ee_k} \vv_h) \Big)
\end{equation}
where $P_h^1$ is the $L^2$ projection on $V^{1,c}_h$.
Our (semi-)discretization of \cref{eqn:weakNS} reads then: 

find $\uu_h \in V^{1,c}_h$ and $p_h \in V^{2,c}_h$ such that
\begin{equation}
    \label{eqn:discrweakNS}
    \int_\Omega \frac{\partial \uu_h}{\partial t} \cdot \vv_h+ \int_\Omega \uu_h \cdot s_h(\uu_h,\vv_h) 
      + \int_\Omega( \Tgrad_h p_h) \cdot \vv_h + \nu \int_\Omega \Tcurl_h \uu_h \ \Tcurl_h \vv_h=0
\end{equation}
holds for all $\vv_h \in V^{1,c}_h$,
together with the incompressibility condition $\Div (\uu_h)=0$ .

\subsection{Pressure equation}
\label{sec:discrpres}
We now describe a way to compute the pressure in order to preserve the divergence condition on $\uu_h$ exactly. 
Assuming a strongly divergence-free velocity at $t=0$, we need
\begin{equation*}
    \frac{\partial }{\partial t} \div(\uu_h)=0
\end{equation*}
which is equivalent to 
\begin{equation*}
    \int_\Omega \frac{\partial }{\partial t} \div(\uu_h) q_h = 0 ~,
    \quad \text{\em i.e.}, \quad
    \int_\Omega \frac{\partial }{\partial t} \uu_h \cdot \Tgrad_h q_h=0
    \qquad ~ \forall q_h \in V^{2,c}_h ~.
\end{equation*}
Using \cref{eqn:discrweakNS} and the relation \eqref{dual_seq_prop},
this leads to defining $p_h$ as the solution to
\begin{equation}
    \label{eqn:discrweakpress}
    \int_\Omega \uu_h \cdot s_h(\uu_h,\Tgrad_h q_h) + \int_\Omega \Tgrad_h p_h \cdot \Tgrad_h q_h=0 
    \qquad  \forall q_h \in V^{2,c}_h ~ .
\end{equation}

Our resulting space discretization for the incompressible Navier-Stokes equations 
is then
\cref{eqn:discrweakNS}--\cref{eqn:discrweakpress}, which we gather here for later purpose: 
find $\uu_h \in V_h^{1,c}$ and $q_h \in V_h^{2,c}$, such that

\begin{subequations}
\label{eqn:confINS}
  \begin{empheq}[left=\empheqlbrace]{align}
    \label{eqn:confINSmom}
    \int_\Omega \frac{\partial \uu_h}{\partial t} \cdot \vv_h+ \int_\Omega \uu_h \cdot s_h(\uu_h,\vv_h) 
        + \int_\Omega (\Tgrad_h p_h) \cdot \vv_h + \nu \int_\Omega \Tcurl_h \uu_h \ \Tcurl_h \vv_h &=0
      \\
    \label{eqn:confINSpres}
    \int_\Omega \uu_h \cdot s_h(\uu_h,\Tgrad_h q_h) + \int_\Omega \Tgrad_h p_h \cdot \Tgrad_h q_h &=0
  \end{empheq}
\end{subequations}
holds for all $\vv_h \in V^{1,c}_h$ and $q_h \in V^{2,c}_h$.

\subsection{Conservation properties}

We now prove that solutions to the semi-discrete problem \cref{eqn:confINS} 
enjoy the same conservation properties 
as the ones shown in \cref{sec:consprop} for the continuous equation.

\begin{proposition}[Preservation of mass]
\label{prop:confpresmass}
Solutions to \cref{eqn:confINS} satisfy

\begin{equation}
\label{eqn:confpresmass}
    \frac{\partial }{\partial t} \Div(\uu_h)=0 ~ .
\end{equation}
Therefore $\Div(\uu_h) = 0 \ \forall t\geq 0$ if $\Div(\uu_h(t=0))=0$ .
\end{proposition}

\begin{proof}[Proof.]
    This follows from the pressure equation computed in the previous section.
\hfill\end{proof}

\begin{proposition}[Preservation of momentum]
\label{prop:confpresmom}
Solutions to \cref{eqn:confINS} preserve the momentum: 
\begin{equation}
    \label{eqn:confmomentconv}
    \frac{\partial}{\partial t}\int_\Omega \uu_h =0 ~ .
\end{equation}
\end{proposition}

\begin{proof}[Proof.]

One can directly adapt the proof of \cref{prop:presmom}, given that we still have the relations $\iTilde^{0h}_{\ee_k} \ee_i = \Tilde{P}_h^0(\delta_{i,k})=\delta_{i,k}$, $\Tgrad_h 1=0$ and $\Tcurl_h \ee_i=0$ by the same argument.
\hfill\end{proof}

\begin{remark}
We here assumed that the constant vector fields are exactly discretized in the space $V^{2,c}_h$. This assumption is satisfied by most of (if not all) the commonly used discretization spaces.
\end{remark}

\begin{proposition}[Preservation of energy]
\label{prop:confpresener}

If $\nu=0$, solutions to \cref{eqn:confINS} preserve the energy: 
\begin{equation}
    \label{eqn:confenergyconv}
    \frac{\partial}{\partial t}\int_\Omega \frac{1}{2}||\uu_h||^2 =0 ~ .
\end{equation}
\end{proposition}

\begin{proof}[Proof.]
The same proof as \cref{prop:presener} is still valid: indeed our discretization preserves the skew-symmetry of $s$, and the pressure term still vanishes by definition of $\Tgrad_h$ as the discrete adjoint of $\Div$ and the fact that $\Div(\uu_h)=0$ due to \cref{prop:confpresmass}.
\hfill\end{proof}

\subsection{Preservation of the velocity invariance property}

If a forcing term is present as in \eqref{eqn:NS_f}, the natural  
approach 
consists of approximating it 
with an $L^2$ projection. Our semi-discrete scheme \eqref{eqn:confINS} then becomes
\begin{subequations}
\label{eqn:confINS_f}
  \begin{empheq}[left=\empheqlbrace]{align}
    \label{eqn:confINSmom_f}
    \int_\Omega \frac{\partial \uu_h}{\partial t} \cdot \vv_h+ \int_\Omega \uu_h \cdot s_h(\uu_h,\vv_h) 
        + \int_\Omega (\Tgrad_h p_h) \cdot \vv_h + \nu \int_\Omega \Tcurl_h \uu_h \ \Tcurl_h \vv_h &= \int_\Omega \ff \cdot \vv_h
      \\
    \label{eqn:confINSpres_f}
    \int_\Omega \uu_h \cdot s_h(\uu_h,\Tgrad_h q_h) 
        + \int_\Omega \Tgrad_h p_h \cdot \Tgrad_h q_h &=
        \int_\Omega \ff \cdot \Tgrad_h q_h
  \end{empheq}
\end{subequations}
for all $\vv_h \in V^{1,c}_h$ and $q_h \in V^{2,c}_h$.
An enjoyable feature of this discretization is that it preserves 
the velocity invariance property \eqref{fip}. 

\begin{proposition}[Discrete velocity invariance property]
\label{prop:fip_h}
The solution to \eqref{eqn:confINS_f} satisfies the following invariance property:
\begin{equation} \label{fip_h}
  \ff \to \ff + \grad \psi  \quad \implies \quad 
  (\uu_h, p_h) \to  (\uu_h, p_h + \Tilde{P}_h^0(\psi))
\end{equation} 
where $\Tilde{P}_h^0$ is the $L^2$ projection on the space $V^{2,c}_h$, see
\eqref{i0h}.
\end{proposition}

\begin{proof}[Proof.]
When the exact source is incremented by $\grad \psi$, the 
right-hand side of the discrete velocity equation \eqref{eqn:confINSmom_f} is incremented by
$$
\int_\Omega \grad \psi \cdot \vv_h 
= -\int_\Omega \psi \div \vv_h 
= -\int_\Omega \Tilde{P}_h^0(\psi) \div \vv_h 
= \int_\Omega \Tgrad_h \Tilde{P}_h^0(\psi) \cdot \vv_h~.
$$
Similarly, the right-hand side of the pressure equation \eqref{eqn:confINSpres_f}
is incremented by $\int_\Omega \Tgrad_h \Tilde{P}_h^0(\psi) \cdot \Tgrad_h q_h$
which allows to prove the claimed property.
\hfill\end{proof}

\subsection{Non-conforming discretization}
\label{sec:nonconf}

We now describe our non-conforming discretization. It relies on broken versions 
of the above spaces, namely spaces $V_h^\ell$ such that 
\begin{equation}
    \label{eqn:nonconfspace}
    V^\ell_h \not \subset V^\ell \quad \text{ for } \quad  \ell=0,1,2~.
\end{equation}
The second ingredient are conforming projection operators 
\begin{equation}
    \label{eqn:confproj}
    P_\ell^c : V^\ell_h \longrightarrow V^\ell_h ~, 
    \qquad (P_\ell^c)^2 = P_\ell^c, \qquad P_\ell^c V^\ell_h = V^{\ell,c}
\end{equation}
which map on conforming discrete spaces defined as 
\begin{equation}
    \label{eqn:confinnonconf}
    V^{\ell,c}_h := V^\ell_h \cap V^\ell~.
\end{equation}

A natural setting is provided by spaces made of local bases  
on disjoint subdomains  
with interface degrees of freedom characterizing the 
conformity conditions, such as classical compatible
spaces of Lagrange/N\'ed\'elec/Raviart-Thomas elements and more generally 
FEEC elements on simplicial meshes \cite{hiptmair2002finite,arnold2018finite}, 
or multipatch spline spaces on mapped cartesian patches \cite{buffa2011isogeometric,gucclu2022broken}.

In this setting, we define 
$$
\Div_h := \Div \circ P^c_1: V^1_h \to V^2_h \quad \text{ and } \quad 
\curl_h := \curl \circ P^c_0: V^0_h \to V^1_h
$$ and $\Tgrad_h: V^2_h \to V^1_h$ (resp. $\Tcurl_h: V^1_h \to V^0_h$) 
is now defined as the adjoint of $-\Div_h$ (resp. $\curl_h$).
We then adapt the conforming scheme \cref{eqn:confINS} by simply replacing 
the conforming (weak) differential operators 
by the non-conforming ones.
Observe that this approach allows to preserve both primal and dual sequence properties, namely
\begin{equation} \label{dual_seq_prop_nc}
  \Div_h \curl_h \vv_h = 0 \qquad \text{ and } \qquad \Tcurl_h \Tgrad_h q_h = 0
\end{equation}
hold for all $\vv_h = 0 \in V^{1}_h$ and $q_h \in V^{2}_h$.
The discrete interior product is then defined as in \eqref{i0h} with 
$\Tilde{P}_h^0$ now denoting the $L^2$ projection on the broken space $V^{2}_h$,
and finally $s_h$ is formally defined as in \eqref{eqn:def_sh}, using the nonconforming
version of the different operators and $P_h^1$ the $L^2$ projection 
on the broken space $V^{1}_h$.
 
This leads to a new scheme on the non-conforming spaces, which takes the same form
as \cref{eqn:confINS},  
namely: find $\uu_h \in V_h^{1}$ and $p_h \in V_h^{2}$ such that

\begin{subequations}
\label{eqn:nonconfINS}
  \begin{empheq}[left=\empheqlbrace]{align}
    \label{eqn:nonconfINSmom}
    \int_\Omega \frac{\partial \uu_h}{\partial t} \cdot \vv_h+ \int_\Omega \uu_h \cdot s_h(\uu_h,\vv_h) + \int_\Omega (\Tgrad_h p_h) \cdot \vv_h + \nu \int_\Omega \Tcurl_h \uu_h \ \Tcurl_h \vv_h &=0 
    \\
    \label{eqn:nonconfINSpres}
    \int_\Omega \uu_h \cdot s_h(\uu_h,\Tgrad_h q_h) + \int_\Omega \Tgrad_h p_h \cdot \Tgrad_h q_h &=0 
  \end{empheq}
\end{subequations}
holds for all $\vv_h \in V^{1}_h$ and $q_h \in V^{2}_h$.

Such a discretization inherits the same conservation properties 
as the conforming one.

\begin{proposition}[Preservation of mass, momentum and energy]
\label{prop:nonconfpres}
Solutions to \cref{eqn:nonconfINS} satisfy
\begin{equation}
\label{eqn:nonconfpresmass}
    \frac{\partial }{\partial t} \Div(\uu_h)=0
\end{equation}
and therefore $\Div(\uu_h)=0 \ \forall t\geq 0$ if $\Div(\uu_h(t=0))=0$,
as well as 
\begin{equation}
    \label{eqn:nonconfmeconv}
    \frac{\partial}{\partial t}\int_\Omega \uu_h =0
    \quad \text{ and } \quad 
    \frac{\partial}{\partial t}\int_\Omega \frac{1}{2}||\uu_h||^2 =0 ~ .
\end{equation}
\end{proposition}

\begin{proof}[Proof.]
For the preservation of mass, one checks that the computations done in \cref{sec:discrpres} are still valid with the addition of the conforming projections $P^i_c$. 
For the momentum, the computations done with $\ee_i$ do not change and 
finally for the energy the skew symmetry of $s_h$ is not affected by the non-conformity
and the pressure term still vanishes due to \cref{eqn:nonconfpresmass}.
\hfill\end{proof}

In the case of a forcing term \eqref{eqn:NS_f}, 
the natural approach following \cite{campos2022broken} 
is to approximate it with an $L^2$ projection on the broken space, 
filtered with the conformin projection $P^c_1$.
Specifically, our semi-discrete scheme \eqref{eqn:nonconfINS} reads in this case
\begin{subequations}
\label{eqn:nonconfINS_f}
  \begin{empheq}[left=\empheqlbrace]{align}
    \label{eqn:nonconfINSmom_f}
    \int_\Omega \frac{\partial \uu_h}{\partial t} \cdot \vv_h+ \int_\Omega \uu_h \cdot s_h(\uu_h,\vv_h) + \int_\Omega (\Tgrad_h p_h) \cdot \vv_h + \nu \int_\Omega \Tcurl_h \uu_h \ \Tcurl_h \vv_h 
      = \int_\Omega \ff \cdot P^1_c \vv_h&
    \\
    \label{eqn:nonconfINSpres_f}
    \int_\Omega \uu_h \cdot s_h(\uu_h,\Tgrad_h q_h) + \int_\Omega \Tgrad_h p_h \cdot \Tgrad_h q_h  = \int_\Omega \ff \cdot P^1_c \Tgrad_h q_h&
  \end{empheq}
\end{subequations}
holds for all $\vv_h \in V^{1}_h$ and $q_h \in V^{2}_h$.

This non-conforming discretization preserves again the velocity invariance property \eqref{fip}. 

\begin{proposition}[Discrete velocity invariance property]
\label{prop:fip_h_nonconf}
The solution to \eqref{eqn:nonconfINS_f} satisfies the following invariance property:
\begin{equation} \label{fip_h_nonconf}
  \ff \to \ff + \grad \psi  \quad \implies \quad 
  (\uu_h, p_h) \to  (\uu_h, p_h + \Tilde{P}_h^0(\psi))
\end{equation} 
where $\Tilde{P}_h^0$ now denotes the $L^2$ projection on the broken space $V^{2}_h$.
\end{proposition}
\begin{proof}[Proof.]
The argument is an extension of the one used for Proposition~\ref{prop:fip_h}.
Here we use that $\Tgrad_h$ is the adjoint of $\div_h = \div P^c_1: V^1_h \to V^2_h$.
When the exact source is incremented by $\grad \psi$, the 
right-hand side of the velocity equation \eqref{eqn:nonconfINSmom_f} 
is incremented by
$$
\int_\Omega \grad \psi \cdot P^c_1 \vv_h 
= -\int_\Omega \psi \div P^c_1 \vv_h 
= -\int_\Omega \Tilde{P}_h^0(\psi) \div_h \vv_h 
= \int_\Omega \Tgrad_h \Tilde{P}_h^0(\psi) \cdot \vv_h~.
$$
Similarly, the right-hand side of the pressure equation \eqref{eqn:nonconfINSpres_f}
is incremented by $\int_\Omega \Tgrad_h \Tilde{P}_h^0(\psi) \cdot \Tgrad_h q_h$
which allows to prove the claimed property.
\hfill\end{proof}

\subsection{Moment preserving conforming projection}

When combined with naive conforming projection operators (for example ones that just average degrees of freedom corresponding to pointwise values on the boundaries), 
the above non-conserving scheme fails to converge with high order. 
A probable cause is that our scheme makes an extensive use of the weak 
operator $\Tgrad_h$ which, as the adjoint of $\Div_h = \Div P^{c}_1$, 
involves the operator 
$P^{c*}_1 : V^1_h \to V^1_h$ defined by 
$\int_\Omega P^{c*}_1 \vv_h \cdot \ww_h := \int_\Omega \vv_h \cdot P^c_1 \ww_h$
for all $\vv_h, \ww_h \in V^1_h$. 
With a basic projection operator, $P^{c*}_1$ is at most of order one, 
which a priori prevents the global scheme to be high order.

One way to recover high order convergence is to use conforming projection operators
$P^c$ that preserve high order polynomial moments, that is: 
\begin{equation} \label{eqn:mpcp}
    \int_\Omega P^c \vv_h \cdot \ww = \int_\Omega \vv_h \cdot \ww 
    \qquad \forall \vv_h \in V^1_h, ~ \ww \in \Po^p. 
\end{equation}

With this property, the adjoint projection $P^{c*}$ is exact on polynomials 
up to degree $\sf{p}$ and is therefore a high order projection.
A simple way to enforce \eqref{eqn:mpcp} is to use an averaging formula with a 
stencil of size $\sim \sf{p}$.
In the case of multipatch spaces with local tensor-product constructions
\cite{gucclu2022broken}, it is enough to describe the 1D construction 
where the conformity amounts to continuity and can be enforced 
by equating interface degrees of freedom corresponding to 
pointwise values on the patch boundaries.
Using broken basis functions $\phi_{k,i}$, $i = 0, \dots, N$ on each patch 
$[x_k,x_{k+1}]$, with boundary values $\phi_{k,i}(x_k) = \delta_{i,0}$
and $\phi_{k,i}(x_{k+1}) = \delta_{i,N}$, a conforming projection with a stencil 
of radius $r \ge 0$ is characterized by its values on the non-conforming basis functions, 
namely $\phi_{k,0}$ and $\phi_{k,N}$. Considering a symmetric construction, we write
\[
P^c_{1d} \phi_{k,0} = \sum_{i=0}^r (c_i \phi_{k,i} + c'_i \phi_{k-1,N-i}), \qquad
P^c_{1d} \phi_{k-1,N} = \sum_{i=0}^r (c'_i \phi_{k,i} + c_i \phi_{k-1,N-i})~.
\]
$P^c_{1d}$ maps on the conforming space if both these functions are continuous, namely
$c_0 = c'_0$, and it is a projection if
$P^c_{1d} (\phi_{k,0} + \phi_{k-1,N}) = \phi_{k,0} + \phi_{k-1,N}$, 
i.e., $c_i + c'_i = \delta_{i,0}$. Thus, 
we must have 
\[
c_0 = c'_0 = \frac 12 \quad \text{and} \quad c_i = -c'_i ~ \text{ for } ~ i > 0~.
\]
A local version of condition \eqref{eqn:mpcp} applied to 
$v_h = \phi_{k,0}$ or $\phi_{k-1,N}$ then yields constraints of the form
\[
\sum_{i=1}^r c_i \int_{[x_k,x_{k+1}]} \phi_{k,i}(x) (x-x_k)^j
= \frac 12 \int_{[x_k,x_{k+1}]} \phi_{k,0}(x) (x-x_k)^j
\quad \text{ for } j = 0, \dots, \sf{p}
\]
which allows us to find coefficients $c_i$ such that \eqref{eqn:mpcp} holds.

\subsection{Adding dissipation through jump penalization}

 In order to stabilize the scheme and limit the non physical jumps artificially created 
 by the use of non conforming spaces we add a symmetric penalization term 
 of the form 
 \newline $\int_\Omega (I-P^c)\uu_h \cdot (I-P^c_1)\vv_h$. 
 The stabilized scheme reads: 

 \begin{subequations}
\label{eqn:nonconfstab}
  \begin{empheq}[left=\empheqlbrace]{align}
    \label{eqn:nonconfstabmom}
    \int_\Omega \frac{\partial \uu_h}{\partial t} \cdot \vv_h + \int_\Omega \uu_h \cdot s_h(\uu_h,\vv_h) + \int_\Omega \Tgrad_h p_h \cdot \vv_h + \nu \int_\Omega \Tcurl_h \uu_h \ \Tcurl_h \vv_h
    \qquad & \\ \nonumber
    + \alpha \int_\Omega (I-P^c_1)\uu_h \cdot (I-P^c_1)\vv_h
    &=0 \\ 
    \label{eqn:nonconfstabpres}
    \int_\Omega \uu_h \cdot s_h(\uu_h,\Tgrad_h q_h) + \int_\Omega \Tgrad_h p_h \cdot \Tgrad_h q_h + \alpha \int_\Omega (I-P^c_1)\uu_h \cdot (I-P^c_1)\Tgrad_h q_h &=0 
  \end{empheq}
\end{subequations}
for all $\vv_h \in V^{1,c}_h$, $q_h \in V^{2,c}_h$.
We observe that the modified pressure equation allows us to preserve the strong
incompressibility condition on the velocity (simply test \cref{eqn:nonconfstabmom} with 
$\vv_h = \Tgrad_h q_h$). If the constant vector fields $\ee_i$ 
belong to the conforming space $V^{1,c}_h$ then the projection property yields 
$(I-P^c_1)\ee_i =0$, hence the momentum is also preserved by the stabilized formulation.
As for the energy, we have the following dissipation rate: 
\begin{equation}
    \frac{\partial}{\partial t}\int_\Omega \frac{|\uu_h|^2}{2} = -\alpha \int_\Omega |(Id-P^c_1)\uu_h|^2~.
\end{equation}
In particular the energy stability is granted, and energy conservation holds for conforming solutions.

\section{Time discretization}
\label{sec:Time_discretization}

\subsection{Crank--Nicolson time scheme}
To preserve energy in the inviscid limit, and dissipate it at a physicaly relevant rate in the viscous case, 
we propose a Crank--Nicolson time discretization.
Starting from an initial solution $\uu_h^0 \in V^{1,c}_h$ satisfying $\Div(\uu_h^0)=0$,
for all $n \ge 1$ we define $\uu_h^{n+1} \in V^{1,c}_h$ and $p_h^{n+1} \in V^{2,c}_h$ 
as the solution to the system
\begin{equation}
\label{eqn:CNINS}
  \left\{\begin{aligned}
      &\int_\Omega \frac{\uu_h^{n+1}-\uu_h^n}{\Delta t} \cdot \vv_h + 
      \int_\Omega \uu_h^{n+\frac 12} \cdot s_h(\uu_h^{n+\frac 12},\vv_h) + \int_\Omega \Tgrad_h p_h^{n+1} \cdot \vv_h 
        \\
      & \mspace{100mu} + \nu \int_\Omega \Tcurl_h \uu_h^{n+\frac 12} \ \Tcurl_h \vv_h 
      + \alpha \int_\Omega (I-P^c_1)\uu_h^{n+\frac 12} \cdot (I-P^c_1)\vv_h =0  \\
    &\int_\Omega \uu_h^{n+\frac 12} \cdot s_h(\uu_h^{n+\frac 12},\Tgrad_h q_h) + \int_\Omega \Tgrad_h p_h^{n+1} \cdot \Tgrad_h q_h = 0
\end{aligned}\right.
\end{equation}
for all $\vv_h \in V^{1,c}_h$ and $q_h \in V^{2,c}_h$,
where we have denoted $\uu_h^{n+\frac 12} := \frac{\uu_h^{n+1}+\uu_h^n}{2}$.

As we will see in Theorem~\ref{thm:cfl} below,
this system is well-posed in a neighborhood of $\uu_h^n$ and
under a CFL condition on the time step we can solve it with a Picard iteration,
see \eqref{eqn:iterCN}.
Moreover, it preserves the previously mentioned invariants. 

\begin{proposition}[Preservation of mass]
\label{prop:CNpresmass}
Solutions to \cref{eqn:CNINS} satisfy
\begin{equation}
\label{eqn:CNpresmass}
    \frac{\Div_h(\uu_h^{n+1})-\Div_h(\uu_h^{n})}{\Delta t}=0
\end{equation}
so that the strong divergence-free property $\Div(\uu_h^n)=0$ holds for all $n$.
    
\end{proposition}

\begin{proof}[Proof.]
We apply \eqref{eqn:CNINS} with $\vv_h = \Tgrad_h q_h$ and use the duality definition of 
$\Tgrad_h$: this yields 
\begin{align*}
    \frac{1}{\Delta t} \big(\int_\Omega \Div_h(\uu_h^{n+1}) - \Div_h(\uu_h^{n}) \big)q_h 
    &= -\int_\Omega \frac{\uu_h^{n+1}-\uu_h^n}{\Delta t} \cdot \Tgrad_h q_h \\
    &= \int_\Omega \uu_h^{n+\frac 12} \cdot s_h(\uu_h^{n+\frac 12},\Tgrad_h q_h) 
      + \int_\Omega \Tgrad_h p_h^{n+1} \cdot \Tgrad_h q_h 
    = 0
\end{align*}
where the second equality follows from $\Tcurl_h \Tgrad_h q_h = 0$, see \eqref{dual_seq_prop} and $\alpha \int_\Omega (I-P^c_1)\uu_h^{n+\frac 12} \cdot (I-P^c_1)\Tgrad_h q_h =0$ because $\Tgrad = P^{c*} \Div^*$.
Since $\Div_h$ maps $V_h^1$ into $V_h^2$, 
this proves \cref{eqn:CNpresmass}.
\hfill\end{proof}

\begin{proposition}[Preservation of momentum]
\label{prop:CNpresmom}
Solutions to \cref{eqn:CNINS} preserves the momentum: 
\begin{equation}
    \label{eqn:CNmomentconv}
    \frac{1}{\Delta t}\Big(\int_\Omega \uu_h^{n+1} - \int_\Omega \uu_h^{n}\Big)=0 ~.
\end{equation}
\end{proposition}

\begin{proof}[Proof.]
We still have $s_h(\uu_h^{n+\frac 12},\ee_l)=\Tgrad_h \iTilde^{0h}_{\ee_l} (\uu_h^{n+\frac 12})$, $\Tcurl \ee_l =0 $ and $ (I-P^c_1)\ee_l =0$, so that we can use the same argument as in the continuous case. Namely, we compute
\begin{align*}
    \int_\Omega \frac{\uu_h^{n+1}-\uu_h^n}{\Delta t} \cdot \ee_l 
    &= - \int_\Omega \uu_h^{n+\frac 12} \cdot s_h(\uu_h^{n+\frac 12},\ee_l) - \int_\Omega \Tgrad_h p_h^{n+1} \cdot \ee_l \\
    &= - \int_\Omega \uu_h^{n+\frac 12} \cdot \Tgrad_h \iTilde^{0h}_{\ee_l} (\uu_h^{n+\frac 12}) 
      - \int_\Omega \Tgrad_h p_h^{n+1} \cdot \ee_l \\
    &= \int_\Omega \Div_h \big(\uu_h^{n+\frac 12}\big) \iTilde^{0h}_{\ee_l} (\uu_h^{n+\frac 12})  
    + \int_\Omega p_h^{n+1} \ \Div_h \ee_l
    = 0
\end{align*}
using that $\uu_h^n$ is strongly divergence free for all $n$.
\hfill\end{proof}

\begin{proposition}[Preservation of energy]
\label{prop:CNpresener}

If $\nu=\alpha=0$, solutions to \cref{eqn:CNINS} preserve the (kinetic) energy: 
\begin{equation}
    \label{eqn:CNenergyconv}
    \frac{1}{\Delta t}\Big(\int_\Omega \frac{1}{2}||\uu_h^{n+1}||^2 - \int_\Omega \frac{1}{2}||\uu_h^{n}||^2\Big) = 0 ~.
\end{equation}
\end{proposition}
\begin{proof}[Proof.]
  Taking $\vv_h = \uu_h^{n+\frac 12}$, we compute
\begin{align*}
    \int_\Omega \frac{1}{2}\frac{||\uu_h^{n+1}||^2-||\uu_h^n||^2}{\Delta t} 
    &= \int_\Omega \frac{\uu_h^{n+1}-\uu_h^n}{\Delta t} \cdot \uu_h^{n+\frac 12}  \\
    &= - \int_\Omega \uu_h^{n+\frac 12} \cdot s_h(\uu_h^{n+\frac 12}, \uu_h^{n+\frac 12}) 
      - \int_\Omega \Tgrad_h p_h^{n+1} \cdot \uu_h^{n+\frac 12}\\
    &= 0 ~,
\end{align*}
using the antisymmetry of $s_h$ and the divergence free condition on $\uu_h^n$ for all time steps.
\hfill\end{proof}

\begin{remark}
In the case of viscous ($\nu > 0$) or broken simulation \eqref{eqn:nonconfstab} with stabilization term $\alpha > 0$, 
the same computation gives us the following rate of dissipation: 
\begin{equation}
    \label{eqn:CNenergydissip}
    \frac{1}{\Delta t}\Big(\int_\Omega \frac{1}{2}||\uu_h^{n+1}||^2 - \int_\Omega \frac{1}{2}||\uu_h^{n}||^2\Big) 
      = - \nu \int_\Omega \big|\Tcurl_h \uu_h^{n+\frac 12}\big|^2 ~ 
        - \alpha \int_\Omega \big|(Id-P^c) \uu_h^{n+\frac 12}\big|^2 ~.
\end{equation}
\end{remark}

\subsection{Iterative procedure for solving the discrete system} 
Since System~\eqref{eqn:CNINS} 
involves a Poisson problem for $p$, we propose to solve it using 
an iterative procedure where $p_h$ and $\uu_h$ are updated successively.

Starting with $\uu_h^{n+1,0}=\uu_h^n$
and writing $\bar \uu_h^{r} := \frac{1}{2}\big(\uu_h^{n+1,r}+\uu_h^n\big)$ 
for $r = 0, 1, \cdots$, we compute
\begin{equation}
\label{eqn:iterCN}
  \left\{ \begin{aligned}
    \int\limits_\Omega (\Tgrad_h p^{n+1,r+1}_h) \cdot (\Tgrad_h q_h) + \int_\Omega \bar \uu_h^{r} \cdot s_h(\bar \uu_h^{r},\Tgrad_h q_h) =0 \\
    \int_\Omega \frac{\uu_h^{n+1,r+1}-\uu_h^n}{\Delta t} \cdot \vv_h+ \int_\Omega \bar \uu_h^{r}\cdot s_h(\bar \uu_h^{r},\vv_h) + \int_\Omega (\Tgrad_h p^{n+1,r+1}_h) \cdot \vv_h \mspace{100mu} & \\
      + \nu \int_\Omega \Tcurl_h \bar \uu_h^{r} \ \Tcurl_h \vv_h 
      + \alpha \int_\Omega (I-P^c_1)\bar \uu_h^{r} \cdot (I-P^c_1)\vv_h &=0
  \end{aligned}\right.
\end{equation}
for all $q_h \in V^{2}_h$ and $\vv_h \in V^{1}_h$, until reaching convergence. 
The converged solution is then solution of \cref{eqn:CNINS}. 
Observe that in this iteration procedure, 
only the pressure step is implicit and corresponds to a Poisson problem.
This problem is solved first and in the absence of Dirichlet boundary conditions the 
well-posedness is ensured by adding a standard regularisation term of the form 
$\ve \int_\Omega p_h q_h$ with very small parameter $\ve$ to ensure that $p_h$
is of zero integral.
The velocity update is explicit, making this iteration procedure computationally efficient. 
One may also note that the divergence free condition and the momentum preservation are achieved at 
every step of this Picard iteration, while energy is preserved only at convergence.

\subsection{Convergence of the Picard iterations and CFL condition}
\label{sec:CFL}

In this section we study the convergence of the iterative scheme \eqref{eqn:iterCN}. 
Let us denote by 
\begin{equation}
  V_h^{div}=\{ \uu \in V_h^{1,c} : \Div \uu =0 \}
\end{equation}
the space of divergence free finite element fields.
Our first observation 
(using the same argument as in \cref{sec:Time_discretization})
is that any solution to System \eqref{eqn:iterCN}
satisfies $\uu_h^{n+1,r+1} \in V_h^{div}$,
and it can be characterized by
\begin{equation}
\label{eqn:iterCN_nopres}
    \int_\Omega \frac{\uu_h^{n+1,r+1}-\uu_h^n}{\Delta t} \cdot \vv_h + \int_\Omega \bar \uu_h^{r} \cdot s_h(\bar \uu_h^{r},\vv_h) + \nu \int_\Omega \Tcurl_h \bar \uu_h^{r} \ \Tcurl_h \vv_h=0 \qquad  \forall \vv_h \in V_h^{div}
\end{equation}
where the pressure gradient term vanishes when tested against divergence free functions. 
The iterative process reads then
$\uu_h^{n+1,r+1} = R_n(\uu_h^{n+1,r})$
where $R_n: V_h^{div} \to V_h^{div}$ is the evolution 
operator defined by 
\begin{equation}
    \label{eqn:evol_op}
    \int_\Omega R_n(\uu_h) \cdot \vv_h = \int_\Omega \uu_h^n \cdot \vv_h - \Delta t \Big( \int_\Omega \frac{\uu_h+\uu_h^n}{2} \cdot s_h\Big(\frac{\uu_h+\uu_h^n}{2},\vv_h\Big) 
      + \nu \int_\Omega \Tcurl_h \frac{\uu_h+\uu_h^n}{2} \ \Tcurl_h \vv_h \Big) 
\end{equation}
for all $\vv_h \in V_h^{div}$. Moreover, 
the solutions to \eqref{eqn:iterCN_nopres} are the fixed points of $R_n$.

We make the following assumption on the discrete spaces, that will be key in the proof.

\begin{assumption}[discrete duality of the $L^\infty$ and $L^1$ norms]
\label{as:dual_L1_Linf} 
    Let $V^*_h = V_h^{div}$, $V^1_h$ or $V^2_h$. There exists 
    a constant $\gamma$ independent of $h$, such that for all $v \in V^*_h$,
    $$
    ||v||_{L^\infty} \leq \gamma \sup_{\substack{w \in V^*_h \\ ||w||_{L^1}\leq 1}} \int_\Omega v \cdot w 
    \qquad \text{and} \qquad     
    ||v||_{L^1} \leq \gamma \sup_{\substack{w \in V^*_h \\ ||w||_{L^\infty}\leq 1}}  \int_\Omega v \cdot w~.
    $$
\end{assumption}

We conjecture that this property holds for regular spline or polynomial finite element spaces on quasi-uniform meshes, using standard scaling arguments and adapting
$L^p$ stability results of $L^2$ projections \cite{crouzeix_stability_1987} 
to handle the subspace of discrete divergence-free functions. 
A detailed proof, however, is not in the scope of the present work. 

Our result then involves the norms
\begin{equation} \label{norms_cd}
||c_h|| := ||c_h||_{L^\infty \times L^\infty \times L^1}
\qquad \text{ and } \qquad ||d_h|| := ||d_h||_{L^\infty \times L^1}  
\end{equation}  
for the operators
$c_h(\uu_h,\vv_h,\ww_h) := \int_\Omega \uu_h \cdot s_h(\vv_h,\ww_h)$ 
and $d_h(\uu_h,\vv_h) := \int_\Omega \Tcurl_h \uu_h \ \Tcurl_h \vv_h$.

\begin{theorem}
  \label{thm:cfl}
  Assume the CFL condition 
  \begin{equation}
      \label{eqn:CFL}
      \frac{\gamma \Delta t}{2} \big(6||c_h||||\uu_h^n||_{L^{\infty}} + \nu ||d_h||\big) < 1 ~.
  \end{equation}
  Then for all $n \ge 1$, the operator $R_n$ defined by \eqref{eqn:evol_op} is a contraction in the $L^\infty$ norm in the domain 
  $U_{\uu_h^n,\delta^n} = \{\uu_h \in V_h^{div} : ||\uu_h-\uu_h^n||_{L^\infty}\leq 2\delta^n ||\uu^n_h||_{L^\infty} \}$ 
  for some $0 < \delta^n \leq 2$.
  In particular, the Picard iterations \eqref{eqn:iterCN} converge towards a solution $\uu^{n+1}_h$
  of System~\eqref{eqn:CNINS}, which   
  satisfies  
  \begin{equation}
  \label{eqn:estimate_iter}
  ||\uu_h^{n+1}-\uu_h^{n}||_{L^{\infty}}\leq 4 ||\uu_h^n||_{L^{\infty}} ~ .
  \end{equation}  
\end{theorem}

We decompose the proof in several lemmas:

\begin{lemma} 
    Under \cref{eqn:CFL} there exists $\delta^n \leq 2$ such that $U_{\uu_h^n,\delta^n}$ is invariant by 
    $R_n$.
\end{lemma}

\begin{proof}[Proof.]

Let $\delta \in \R$ and $\uu_h \in U_{\uu_h^n,\delta}$.
For any $\vv_h \in V_h^{div}$, we have 
\begin{align}
    \int_\Omega \big(R_n(\uu_h)-\uu_h^n\big) \cdot \vv_h 
    &= \Delta t \Big(c_h\Big(\frac{\uu_h+\uu_h^n}{2},\frac{\uu_h+\uu_h^n}{2},\vv_h\Big) + \nu d_h\Big(\frac{\uu_h+\uu_h^n}{2},\vv_h\Big) \Big) \\
    &\leq \Delta t ||\vv_h||_{L^1} \Big( ||c_h|| \Big\|\frac{\uu_h+\uu_h^n}{2}\Big\|_{L^{\infty}}^2+\nu ||d_h|| \Big\|\frac{\uu_h+\uu_h^n}{2}\Big\|_{L^{\infty}} \Big) ~.
\end{align}
By Assumption~\ref{as:dual_L1_Linf} this yields 
\begin{align}
    ||R_n(\uu_h)-\uu_h^n||_{L^{\infty}} 
    &\leq \gamma \sup_{||\vv_h||_{L^1}\leq 1} \int_\Omega \big(R_n(\uu_h)-\uu_h^n \big) \cdot \vv_h \\
    &\leq \gamma \Delta t \Big( ||c_h||\Big\|\frac{\uu_h+\uu_h^n}{2}\Big\|_{L^{\infty}}^2+\nu ||d_h||\Big\|\frac{\uu_h+\uu_h^n}{2}\Big\|_{L^{\infty}}\Big) \\
    &\leq \gamma \Delta t ||\uu_h^n||_{L^{\infty}} \big( ||c_h||(1+\delta)^2||\uu_h^n||_{L^{\infty}}+\nu ||d_h||(1+\delta)\big) ~ . 
\end{align}
For $R_n(\uu_h)$ to be in $U_{\uu_h^n,\delta}$, one thus needs 
$\gamma \Delta t \big( ||c_h||(1+\delta)^2||\uu_h^n||_{L^{\infty}}+\nu ||d_h||(1+\delta)\big) \leq 2\delta$.
Let $\lambda = \gamma \Delta t$, $\alpha = ||c_h||||\uu_h^n||_{L^{\infty}}$, $\beta=\nu ||d_h||$. The condition is: 
$\lambda \big( \alpha (1+\delta)^2+\beta (1+\delta)\big) - 2\delta \leq 0$, equivalent to: $\lambda \big( \alpha \delta^2 +(2 \alpha+ \beta -\frac{2}{\lambda})\delta +\alpha+\beta \big)  \leq 0$.

We therefore need to study the polynomial $Q(X)=\alpha X^2 +(2 \alpha+ \beta -\frac{2}{\lambda})X +\alpha+\beta$. Recall that our aim is to find $\delta>0$ such that $Q(\delta) \le 0$.

First, we infer from $\alpha, \beta >0$ that the coefficient of $X^2$ is positive, so that for such a $\delta$ to exist $Q$ needs to have at least one positive root. Second, from $Q(0)>0$ we infer that (if $Q$ has roots) both roots have the same sign, which is the sign of the argminimum of $Q$, namely $\frac{1}{\alpha}(\frac{2}{\lambda}-2 \alpha -\beta)$:
since $\lambda, \alpha>0$, a first condition is thus

\begin{equation}
\label{eqn:firstCFL}
    \lambda(\alpha+\frac{\beta}{2})<1 ~.
\end{equation}

This first condition is clearly weaker then the claimed CFL \cref{eqn:CFL}.
We are then left to verify that Q has roots, and for this we need the discriminant 
$\Delta = (2 \alpha+ \beta -\frac{2}{\lambda})^2-4\alpha(\alpha+\beta)$ to be positive. 
We have $\Delta = (2 \alpha+ \beta -\frac{2}{\lambda}-2\sqrt{\alpha(\alpha+\beta)})(2 \alpha+ \beta -\frac{2}{\lambda}+2\sqrt{\alpha(\alpha+\beta)})$.
The first term is clearly negative (due to the previous condition) so the condition $\Delta>0$ is equivalent to $2 \alpha+ \beta -\frac{2}{\lambda}+2\sqrt{\alpha(\alpha+\beta)}<0$, i.e., $ \lambda(\alpha+ \frac{\beta}{2} +\sqrt{\alpha(\alpha+\beta)})<1$. Since $\sqrt{\alpha(\alpha+\beta)}\leq \alpha + \frac{\beta}{2}$ we have the stronger condition:
\begin{equation}
\label{eqn:secondCFL}
    \lambda(\alpha+ \frac{\beta}{2})<\frac{1}{2}~.
\end{equation}
Under this condition (which is still weaker than \cref{eqn:CFL}), 
there exists $\delta^n > 0$ such that $U_{\uu_h^n,\delta^n}$ is stable by $R_n$.
We choose for $\delta^n$ the smallest root of $Q$.
Let us estimate its value:
for this we rewrite $Q(X)=aX^2+bX+c$, 
and from the previous computations we know that $a>0, b<0, c>0, \Delta = b^2-4ac>0$ 
so that $0 \leq \frac{4ac}{b^2} \leq 1$. We have 
$$ 
\delta^n = \frac{-b-\sqrt{b^2-4ac}}{2a}=\frac{-b}{2a}\Big(1-\sqrt{1-\frac{4ac}{b^2}}\Big)\leq \frac{-b}{2a}\frac{4ac}{b^2}=\frac{-2c}{b}
$$
(where we have used the fact that for $0\leq x \leq 1$, $x \geq 1-\sqrt{1-x}$).
Replacing $a$, $b$ and $c$ gives
$\delta^n \leq 
\lambda \frac{\alpha+\beta}{1-\lambda(\alpha+\frac{\beta}{2})}$.
Now, $\lambda(\alpha+ \frac{\beta}{2})<\frac{1}{2}$, so that $1-\lambda(\alpha+ \frac{\beta}{2})>\frac{1}{2}$.
This yields 
\begin{equation} \label{bound_delta}
\delta^n \leq \lambda \frac{\alpha+\beta}{1-\lambda(\alpha+\frac{\beta}{2})}\leq 2 \lambda(\alpha+\beta)\leq 2 \lambda(2\alpha+\beta)\leq 2~.  
\end{equation}
\hfill\end{proof}

\begin{lemma}
    Assume the CFL condition \cref{eqn:CFL}. Then the operator $R_n$ is contracting in 
    $U_{\uu_h^n,\delta^n}$ for $\delta^n$ defined as the smallest root in the proof of the previous lemma.
\end{lemma}

\begin{proof}[Proof.]
    
Let $\uu_h, \uu_h' \in U_{\uu_h^n,\delta^n}$, then
$||R_n(\uu_h)-R_n(\uu'_h)||_{L^{\infty}} \leq \gamma \sup_{||\vv_h||_{L^1}\leq 1} \int_\Omega \big(R_n(\uu_h)-R_n(\uu'_h)\big) \cdot \vv_h$.
We observe that by tri-linearity of $c_h$, one as
$$
c_h\Big(\frac{\uu_h+\uu_h^n}{2},\frac{\uu_h+\uu_h^n}{2},\vv_h\Big)
  -c_h\Big(\frac{\uu'_h+\uu_h^n}{2},\frac{\uu'_h+\uu_h^n}{2},\vv_h\Big)
  = c_h\Big(\frac{\bar \uu_h+\uu_h^n}{2},\tilde \uu_h,\vv_h\Big)
    +c_h\Big(\tilde \uu_h,\frac{\bar \uu_h+\uu_h^n}{2},\vv_h\Big)
    $$
where we have denoted $\bar \uu_h = \frac{\uu_h+\uu_h'}{2}$ and $\tilde \uu_h = \frac{\uu_h-\uu_h'}{2}$.
Then, for any $\vv_h \in V_h^{div}$, we have 
\begin{align}
    \int_\Omega \big(R_n(\uu_h)-R_n(\uu'_h)\big) \cdot \vv_h 
      &= \Delta t \Big(c_h\Big(\frac{\bar \uu_h+\uu_h^n}{2},\tilde \uu_h,\vv_h\Big)
        +c_h\Big(\tilde \uu_h,\frac{\bar \uu_h+\uu_h^n}{2},\vv_h\Big) +\nu d_h(\tilde \uu_h,\vv_h) \Big) \\
    & \leq  \Delta t \big(||c_h||||\bar \uu_h+\uu_h^n||_{L^{\infty}}||\tilde \uu_h||_{L^{\infty}}||\vv_h||_{L^1} + \nu||d_h||||\tilde \uu_h||_{L^{\infty}}||\vv_h||_{L^1}\big) \\
    & \leq  \Delta t \big(||c_h||(2+2\delta^n)||\uu_h^n||_{L^{\infty}} + \nu||d_h||\big)
      \frac{||\uu_h-\uu_h'||_{L^{\infty}}}{2}||\vv_h||_{L^1}~,
\end{align}
so that using \eqref{bound_delta} yields
\begin{equation}
    ||R_n(\uu_h)-R_n(\uu'_h)||_{L^{\infty}}\leq \frac{\gamma\Delta t}{2} \big(6 ||c_h|| ||\uu_h^n||_{L^{\infty}} + \nu ||d_h||\big)||\uu_h-\uu_h'||_{L^{\infty}}
\end{equation}
and $R_n$ is indeed contracting on $U_{\uu^n,\delta^n}$ due to \cref{eqn:CFL}.
\hfill\end{proof}

We can now conplete the argument.
\begin{proof}[Proof of Theorem~\ref{thm:cfl}.]
Now that the contraction property of $R_n$ is established,
the convergence of the iterative procedure follows from the Banach fixed-point theorem. 
The estimate \eqref{eqn:estimate_iter} is a reformulation that the limit must be in $U_{\uu_h^n,\delta^n}$.

\hfill\end{proof}

\subsection{Scaling of the CFL condition}
\label{sec:scale_CFL}

In this section  
we study how the CFL condition scales with the discretization parameters. 
For this we begin by meshing the domain $\Omega$ with disjoint elements, namely
\begin{equation}
  \label{eqn:mesh_assumption}
  \Omega=\bigcup\limits_{i=1}^N K_i ~ \quad \text{ with } \quad 
  |K_i \cap K_j|=0 ~ \quad \forall i \neq j~. 
\end{equation}
We also assume that there is a reference element $\hK$ and mappings $(\phi_i)_{i=1,...,N}$  such that for all $i$, $K_i=\phi_i(\hK)$. 
For simplicity we assume that the mappings are affine, 
\begin{equation}
\label{eqn:form_phi}
\phi_i(x,y)=(h_x^i x +b_x^i, h_y^i y + b_y^i) ~ .
\end{equation}
We denote $h=\max_{i=1,...,N} \max (h^i_x,h^i_y)$ and assume that the mesh is regular in the sense that there is a constant $\sigma$ (independent of $h$) such that 
$h_x^i, h_y^i \geq \sigma h$ for all $i$. 
For later purpose we observe that the Jacobian matrix of $\phi_i$ reads 
    \begin{equation}
        \label{eqn:dphi}
        D\phi_i= 
        \begin{pmatrix}
        h_x^i & 0 \\
        0 & h_y^i
        \end{pmatrix} ~ ,
    \end{equation}
    so that for all vector $\uu \in \R^2$, it holds
    \begin{equation}
        \label{eqn:estim_Dphiu}
        \sigma h |\uu| \leq |D \phi_i \uu| \leq h |\uu| ~ .
    \end{equation}
    Moreover the Jacobian determinant $J_{\phi_i} = h_x^i h_y^i$ satisfies
    \begin{equation}
        \label{eqn:bounds_Jac}
        \sigma^2 h^2 \leq J_{\phi_i} \leq h^2 ~ .
    \end{equation}

We next suppose that we are given a reference de Rham sequence 
\begin{equation}
\label{eqn:refDeRham}
\xymatrix{
 \hV^0\subset V^0(\hK) \ar[r]^-{\curl} & \hV^1\subset V^1(\hK) \ar[r]^-{div} & \hV^2\subset V^2(\hK)
 }~.
\end{equation}
According to \cite{hiptmair2002finite,buffa2011isogeometric}, 
we can define local de Rham sequences on the elements $K_i$
by transporting the references spaces \cref{eqn:refDeRham} 
with the push-forward operators 
\begin{align}
    \label{eqn:defF_i0}
    F_i^0(\huu) &=\huu\circ \phi_i^{-1} ~ ,\\
    \label{eqn:defF_i1}
    F_i^1(\huu) &=(J_{\phi_i}^{-1}D\phi_i \huu)\circ \phi_i^{-1} ~ ,\\
    \label{eqn:defF_i2}
    F_i^2(\huu) &=(J_{\phi_i}^{-1}\huu)\circ \phi_i^{-1} ~ .
\end{align}
The resulting local spaces are of the form
\begin{equation}
    \label{eqn:physical_spaces}
    V_h^k(K_i)=\{F_i^k(\huu), \uu \in \hV^k \} ~ ,
\end{equation}
so that we can finaly define our global broken spaces as 
\begin{equation}
    \label{eqn:full_space}
    V_h^k = \{v \in V^k(\Omega) : \forall i, \ v_{|K_i}\in V_h^k(K_i)\} ~.
\end{equation}
Because the divergence operator commutes with the $F_i^1$ push-forward operators,
we also have 
\begin{equation}
    \label{eqn:global_Vdiv}
    V_h^{div} = \{\uu \in V^1(\Omega) : \forall i, \ \uu_{|K_i}\in F_i(\hV^{div})\} ~,
\end{equation}
where we have defined $\hV^{div} = \{ \huu \in \hV : \Div(\huu)=0 \}$.

We can now study the scaling of $||c_h||$ and $||d_h||$ with $h$:

\begin{proposition}
\label{prop:estim_c_d}
    If the conforming projectors $P_c^k$ are uniformly bounded (with respect to $h$) in the $L^1$ and $L^\infty$ norms, then there exists constants $C$ and $D$ independent of $h$ such that 
    \begin{equation} \label{eqn:scaling_c_d}
      ||c_h||\leq \frac{C}{h}
      \qquad \text{ and } \qquad
      ||d_h||\leq \frac{D}{h^2} ~.
    \end{equation} 
\end{proposition}
\begin{proof}[Proof.]
    We must bound $c_h(\uu_h,\vv_h,\ww_h)= \frac{1}{2}\int_\Omega \uu_h \cdot \big( \sum\limits_{k=1}^{d} (\iTilde^{0h}_{\ee_k} \ww_h)\Tgrad_h (\iTilde^{0h}_{\ee_k} \vv_h) - (\iTilde^{0h}_{\ee_k} \vv_h)\Tgrad_h (\iTilde^{0h}_{\ee_k} \ww_h) \big)$
    for $\uu_h,\vv_h \in L^\infty$ and $\ww_h \in L^1$. 
    We begin with the second term and bound  
    \begin{equation} \label{boundc1}
      \int_\Omega \uu_h \cdot (\iTilde^{0h}_{\ee_k} \vv_h) \Tgrad_h( \iTilde^{0h}_{\ee_k} \ww_h) \leq ||\uu_h||_{L^\infty} ||\iTilde^{0h}_{\ee_k} \vv_h||_{L^\infty} ||\Tgrad_h \iTilde^{0h}_{\ee_k}\ww_h||_{L^1}~.
    \end{equation} 
    By definition of the discrete interior product \eqref{i0h} 
    (and its extension in the non-conforming case) we have for 
    all $w_h \in V^2_h$
    \begin{equation*}
        \int_\Omega \iTilde^{0h}_{\ee_k} \vv_h w_h 
        = \int_\Omega (\vv_h)_k w_h 
        \leq ||(\vv_h)_k||_{L^\infty} ||w_h||_{L^1}
        \leq ||\vv_h||_{L^\infty} ||w_h||_{L^1} ~,
    \end{equation*}
    hence Assumption~\ref{as:dual_L1_Linf} yields  
    $||\iTilde^{0h}_{\ee_k} \vv_h||_{L^\infty}\leq \gamma ||\vv_h||_{L^\infty}$.

In a similar fashion, we can use 
$||\Tgrad_h \iTilde^{0h}_{\ee_k}\ww_h||_{L^1}\leq c \sup_{||\vv'_h||_{L^\infty}=1} \int_\Omega \vv'_h \cdot \Tgrad_h \iTilde^{0h}_{\ee_k}\ww_h$. 
For $\vv'_h \in V^1_h$ with $||\vv'_h||_{L^\infty}=1$, we write 
$\vv_h=P_c^1 \vv'_h$
and further denote the pull-backs on a cell $K_i$ by 
$\hvv^i_h := (F^1_i)^{-1}(\vv_h)$ and $\hww^i_h := (F^1_i)^{-1}(\ww_h)$.
We compute:
\[
\begin{aligned}
    \int_\Omega \vv'_h \cdot \Tgrad_h \iTilde^{0h}_{\ee_k}\ww_h
    &= -\int_\Omega (\Div \vv_h) \iTilde^{0h}_{\ee_k}\ww_h
    = -\int_\Omega (\Div \vv_h) (\ww_h)_k
    = -\sum\limits_{i=1}^N \int_{K_i} (\Div \vv_h) (\ww_h)_k \\
    &= -\sum\limits_{i=1}^N \int_{K_i} (\Div F^1_i(\hvv^i_h)) (F^1_i(\hww^i_h))_k
    = -\sum\limits_{i=1}^N \int_{K_i} (F^2_i(\Div \hvv^i_h)) (F^1_i(\hww^i_h))_k
    \\
    &= - \sum\limits_{i=1}^N \int_{\hat K} (\Div \hvv^i_h) (J_{\phi_i}^{-1} D_{\phi_i} \hww^i_h)_k 
    \leq \sigma^{-2} h^{-1} \norm{\widehat\Div}_{L^\infty} \sum\limits_{i=1}^N \norm{\hvv^i_h}_{L^\infty} \norm{(\hww^i_h)_k}_{L^1}~. 
\end{aligned}
\]
Here we have used the bounds \eqref{eqn:estim_Dphiu}--\eqref{eqn:bounds_Jac} and $\norm{\widehat\Div}_{L^\infty}$ is the norm of the $\Div$ operator from 
$(\hV^1,L^\infty)$ to $(\hV^2,L^\infty)$, which is independent of $h$.
Using the form of the push-forward \cref{eqn:defF_i0,eqn:defF_i1,eqn:defF_i2}
and again the bounds \eqref{eqn:estim_Dphiu}--\eqref{eqn:bounds_Jac} we find 
$$
||\hvv^i_h||_{L^\infty} \le \sigma^{-1} h ||\vv_h||_{L^\infty}
\quad \text{ and } \quad 
 \sum\limits_{i=1}^N ||\hww^i_h||_{L^1} \le 
 (\sigma h)^{-1} \sum\limits_{i=1}^N ||\ww_h||_{L^1(K_i)}  \le (\sigma h)^{-1}  ||\ww_h||_{L^1}~.
$$
It follows that
\begin{equation}
\int_\Omega \vv'_h \cdot \Tgrad_h \iTilde^{0h}_{\ee_k}\ww_h 
  \leq (\sigma^{4} h)^{-1} ||\vv_h||_{L^\infty} ||\ww_h||_{L^1} 
  \le (\sigma^{4} h)^{-1} ||P^1_c ||_{L^\infty} ||\vv'_h||_{L^\infty} ||\ww_h||_{L^1}   
\end{equation}
where $\norm{P_c^1}_{L^\infty}$ is the norm of
the conforming projection:
as the latter consists of local averages of local basis functions, 
it is uniformly bounded with respect to $h$. 
Using $||\vv'_h||_{L^\infty} = 1$, we find
\begin{equation}
||\Tgrad_h \iTilde^{0h}_{\ee_k}\ww_h||_{L^1}\leq c h^{-1} ||\ww_h||_{L^1}
\end{equation}
with a constant independent of $h$. Plugging these estimates in \eqref{boundc1} gives
\begin{equation}
    \int_\Omega \uu_h \cdot (\iTilde^{0h}_{\ee_k} \vv_h)\Tgrad_h (\iTilde^{0h}_{\ee_k} \ww_h) \leq c h^{-1} ||\uu_h||_{L^\infty} ||\vv_h||_{L^\infty} ||\ww_h||_{L^1}~.
\end{equation}
Repeating the same computations and switching the $L^1$ and $L^\infty$ norms gives
\begin{equation}
    \int_\Omega \uu_h \cdot \iTilde^{0h}_{\ee_k} \ww_h \Tgrad_h \iTilde^{0h}_{\ee_k} \vv_h \leq c h^{-1}||\uu_h||_{L^\infty} ||\vv_h||_{L^\infty} ||\ww_h||_{L^1}
\end{equation}
with another constant independent of $h$,
hence the bound for $||c_h||$ in \eqref{eqn:scaling_c_d}.
The same arguments allow to prove that $||\Tcurl \vv_h||_{L^\infty} \leq \frac{c}{h}||\vv_h||_{L^\infty}$ and $||\Tcurl \vv_h||_{L^1} \leq \frac{c}{h}||\vv_h||_{L^1}$, 
which justifies the bound for $||d_h||$ in \eqref{eqn:scaling_c_d}.
\hfill\end{proof}

As a consequence we find that a sufficient condition for  
the CFL inequality \eqref{eqn:CFL} takes the form
\begin{equation}
    \label{eqn:CFL_scaling}    
    C' \Delta t \Big(\frac{||\uu_h^n||_{L^{\infty}}}{h} + \frac{\nu}{h^2}\Big) < 1
\end{equation}
with a constant $C'$ independent of $h$. 
    
\section{Adding boundary conditions}
\label{sec:BC}
We now focus on integrating boundary conditions to our numerical scheme,
in order to solve problems of the form  
\begin{equation}
    \label{eqn:NS_bound}
    \left\{
    \begin{array}{rll}
        \frac{\partial \uu}{\partial t}+ (\uu \cdot \grad) \uu + \grad p - \nu \Delta \uu & = 0 & \text{on } \Omega   \\
         \div \uu & = 0 & \text{on } \Omega   \\
         \uu \times \nn & = u_{t} & \text{on } \Gamma_t \\
         \uu \cdot \nn & = u_{n} & \text{on } \Gamma_n \\
         p & = p_b & \text{on } \Gamma_p
    \end{array}
    \right.
\end{equation}
where $\Gamma_t$, $\Gamma_n$ and $\Gamma_p$ are three pieces of 
the boundary $\partial \Omega$ such that $\Gamma_n$ and $\Gamma_p$ 
are disjoint and cover $\partial \Omega$.  

The implementation of boundary conditions is done in two differents ways: 
normal boundary conditions on $\uu \cdot \nn$ are imposed in a \emph{strong} form,
since we use a discrete $H(\div)$ space for the velocity where boundary degrees of 
freedom correspond to discrete fluxes.  
On the other hand, boundary conditions on $p$ and on the tangential component
$\uu \times \nn$ are integrated in \emph{weak} form, that is replacing term where they appear on boundary integrals. 
Since the tangential velocity on a boundary is often related to the viscosity, it is natural to make it appear on the vector Laplacian, hence 
in a weak definition of the discrete curl operator. The pressure condition is in some sense a dual condition to the normal velocity condition 
(which justifies that $\Gamma_n$ and $\Gamma_p$ form a partition of the boundary) 
therefore it is naturally treated in a weak form. 

For simplicity, in this section 
we will not separe the description of conforming and non-conforming case as they can be treated in the same way. Thus, we assume here $V^\ell_h = V^{\ell,c}_h$.

\subsection{Discrete differential operators with boundary conditions}
A crucial point to understand the implementation of weak boundary conditions is the mutual relation of primal/dual differential operators via integration by parts. Every time a differential operator is integrated by parts, a corresponding boundary integral must be considered. 
In particular, in the absence of periodic or homogeneous boundary conditions
the weak differential operators must be redefined:
Specifically, the discrete operators 
$$
\Tgrad_h: V^{2}_h \to V^{1}_h \qquad \text{ and } \qquad
\Tcurl_h: V^{1}_h \to V^{0}_h
$$
are defined in full generality by the relations
\begin{equation} \label{eqn:dualh_eq_fg}
  \left\{ \begin{aligned}
  \int_\Omega \vv_h \cdot \Tgrad_h q_h  &= - \int_\Omega \div(\vv_h) q_h 
      + \int_{\partial \Omega} q_h (\vv_h \cdot \nn)
  \\
  \int_\Omega (\Tcurl_h \vv_h ) \omega_h &= \int_\Omega \vv_h \cdot \curl\omega_h 
      - \int_{\partial \Omega} (\vv_h \times \nn)\omega_h
  \end{aligned} \right.
\end{equation}
for all $q_h \in V^{2}_h$, $\vv_h \in V^{1}_h$, and all $\omega_h \in V^0_h$.

As discussed above, normal velocity boundary conditions will be handled in strong
form. For this we introduce the space $\mathring{V_h^1}$ of velocity fields
with \emph{homogenous} fluxes on $\Gamma_{n}$, and suppose that we are given a projection $P_n$ on it,  i.e.
$$
\mathring{V_h^1}=\{\vv_h \in V_h^1 | \vv_h \cdot \nn =0 \text{ on } \Gamma_{n}\}, \quad\quad P_n : V_h^1 \xrightarrow{} \mathring{V_h^1}.
$$ 
Specifically, $P_n$ cancels flux degrees of freedom on the boundary $\Gamma_n$
and does not modify the other ones, so that we have 
\begin{equation} \label{Pn_on_Gp}
P_n\vv_h \cdot \nn = \vv_h \cdot \nn  
\quad  \text{ on } \quad \Gamma_p = \partial \Omega \setminus \Gamma_n.
\end{equation}
To account for the pressure and tangential velocity boundary conditions in \eqref{eqn:NS_bound}, we next define two discrete operators
$\Tgrad_{h,p_b}: V^{2}_h \to V^{1}_h$ and 
$\Tcurl_{h,u_t}: V^{1}_h \to V^{0}_h$ which involve the boundary terms 
$p_b$ and $u_t$. 
Because in our velocity equation the discrete gradients will be tested
against functions of the form $P^*_n \vv_h$, we consider
an operator $\Tgrad_{h,p_b}$ that actually maps into the subspace $P^*_n V^1_h$.
We define 
\begin{equation} \label{eqn:dualh_eq_bc}
  \left\{ \begin{aligned}
  \int_\Omega \vv_h \cdot \Tgrad_{h,p_b} q_h  &= - \int_\Omega \div(P_n \vv_h) q_h 
      + \int_{\Gamma_p} p_b (\vv_h \cdot \nn)
  \\
  \int_\Omega (\Tcurl_{h,u_t} \vv_h )\omega_h &= \int_\Omega \vv_h \cdot \curl\omega_h 
      - \int_{\partial \Omega \setminus \Gamma_t} (\vv_h \times \nn) \omega_h
      - \int_{\Gamma_t} u_t \omega_h~.
  \end{aligned} \right.
\end{equation}
We observe that this definition indeed yields
$\Tgrad_{h,p_b}: V^{2}_h \to P^*_n V^{1}_h$, since 
$$
\int_\Omega \vv_h \cdot (I-P^*_n) \Tgrad_{h,p_b} q_h = 
\int_\Omega (I-P_n)\vv_h \cdot \Tgrad_{h,p_b} q_h = 
\int_{\Gamma_p} p_b \big(((I-P_n)\vv_h) \cdot \nn\big) = 0
$$ 
for all $\vv_h \in V^1_h$, using \eqref{Pn_on_Gp}.
It will also be convenient to denote by 
$$
\Bgrad_h: V^{2}_h \to V^{1}_h \qquad \text{ and } \qquad
\Bcurl_h: V^{1}_h \to V^{0}_h
$$
the \lq\lq boundaryless\rq\rq\ $L^2$ adjoints of the primal div and curl operators 
\begin{equation} \label{eqn:dualh_eq_bl}
  \left\{ \begin{aligned}
  \int_\Omega \vv_h \cdot \Bgrad_h q_h
    &:= -\int_\Omega \div(\vv_h) q_h
    = \int_\Omega \vv_h \cdot \Tgrad_h q_h - \int_{\partial \Omega} q_h (\vv_h \cdot \nn) 
  \\
  \int_\Omega (\Bcurl_h \vv_h )\omega_h 
  &:= \int_\Omega \vv_h \cdot \curl\omega_h 
  = \int_\Omega (\Tcurl_h \vv_h) \omega_h 
    + \int_{\partial \Omega} (\vv_h \times \nn) \omega_h
  \end{aligned} \right.
\end{equation}

\subsection{Discrete velocity equation with boundary conditions}

In the absence of periodic or homogenous boundary conditions, 
the discrete operator $s_h$ defined by \eqref{eqn:def_sh} should 
involve the discrete gradient operator {\em with boundary terms} 
defined in the previous section. 
Moreover, the integration by parts 
leading to the bilinear operator $s$ at the continuous level also 
involves boundary terms: it reads 
\begin{equation*} 
    \int_\Omega \uu \cdot s(\uu, \vv) = \int_\Omega(\uu \cdot \grad) \uu \cdot \vv = \frac{1}{2} \big(
    \int_\Omega(\uu \cdot \grad) \uu \cdot \vv -\int_\Omega(\uu \cdot \grad) \vv \cdot \uu + \int_{\partial \Omega} (\uu \cdot \vv ) (\uu \cdot \nn) \big).
\end{equation*}
As a consequence, using the operators from \eqref{eqn:dualh_eq_fg} and \eqref{eqn:dualh_eq_bl} we define 
\begin{equation} 
  s_h(\uu_h,\vv_h) := 
  \frac{1}{2} \big(\sum\limits_{k=1}^{d} \iTilde^{0h}_{\ee_k} \vv_h \Tgrad_h \iTilde^{0h}_{\ee_k} \uu_h - \iTilde^{0h}_{\ee_k} \uu_h \Bgrad_h \iTilde^{0h}_{\ee_k} \vv_h \big).
\end{equation}

To handle the viscosity term we rewrite as above
$-\Delta \uu = \curl \omega$ with $\omega = \curl \uu$. 
For a discrete velocity $\uu_h \in V^1_h$ this is now consistently approximated as 
$-\Delta_h \uu_h \in V^1_h$, defined by 
\begin{align}\label{eq:Laplacebnd2}
\left\{ \begin{aligned}
\omega_h &:= \Tcurl_{h,u_t} \uu_h  \in V^0_h \\  
-\int_\Omega \Delta_h \uu_h \cdot \vv_h &= \int_\Omega (\curl \omega_h) \cdot \vv_h = \int_\Omega \Tcurl_{h,u_t} \uu_h  \Bcurl \vv_h
	\end{aligned}\right.
\end{align}
for all $\vv_h \in V^1_h$, using the operators from \eqref{eqn:dualh_eq_bc} and \eqref{eqn:dualh_eq_bl}.

To finally derive an equation for the discrete velocity 
we decompose the solution in the form 
$
\uu_h = \buu_h + \uue_h
$ with
$\uue_h := P_n \uu_h \in \mathring{V_h^1}$. 
In this decomposition, $\buu_h$ 
carries the $\uu \cdot \nn$ boundary condition and is supposed to be constant in time. 
By doing this, we have restricted the space for our unknown to $\mathring{V_h^1}$ and, therefore, we only need to test our equation against functions of the form 
$P_n^* \vv_h$  
with $\vv_h \in V_h^1$. 
Using the discrete operators defined above, this gives
\begin{multline} \label{eq:uue_t}
    \int_\Omega \frac{\partial \uue_h}{\partial t} \cdot P_n^* \vv_h+ \int_\Omega (\buu_h + \uue_h) \cdot s_h(\buu_h + \uue_h,P_n^*\vv_h) + \int_\Omega \Tgrad_{h,p_b} p_h \cdot P_n^* \vv_h 
    \\
    + \nu \int_\Omega \Tcurl_{h,u_t} (\buu_h + \uue_h) \cdot \Bcurl_h P_n^* \vv_h=0 
\end{multline}
for all $\vv_h \in V^1_h$. Moreover, using $(I-P_n) \buu_h 
= (I-P_n)^2 \uu_h = \buu_h$ our condition $\frac{\partial \buu_h}{\partial t} = 0$ yields
\begin{equation} \label{eq:buu_t}
  \int_\Omega \frac{\partial \buu_h}{\partial t} \cdot (I-P_n^*) \vv_h = 
  \int_\Omega \frac{\partial \buu_h}{\partial t} \cdot \vv_h = 0.  
\end{equation}
Adding \eqref{eq:uue_t} and \eqref{eq:buu_t} leads to the discrete velocity equation
\begin{equation} \label{uh_eq_bc}
    \int_\Omega \frac{\partial \uu_h}{\partial t} \cdot \vv_h
      + \int_\Omega \uu_h \cdot s_h(\uu_h,P_n^*\vv_h) 
      + \int_\Omega \Tgrad_{h,p_b} p_h \cdot P_n^* \vv_h 
    + \nu \int_\Omega (\Tcurl_{h,u_t} \uu_h) (\Bcurl_h P_n^* \vv_h)=0 
\end{equation}
In particular we obtain an equation written in the full space $V_h^1$, 
which avoids implementing the homogeneous condition within the finite element space. 
This is also beneficial in the case of inhomogeneous boundary conditions 
(as will be considered in \cref{sec:Numerics}), 
where the solution is searched for in the full space.

\subsection{Discrete pressure equation with boundary conditions}
\label{sec:pres_eq_bnd}

For the pressure, we use the same strategy as before (\cref{sec:discrpres}), that is, finding a pressure equation that leads to divergence preservation. 
Since 
$\Bgrad_h$ has been defined as the discrete adjoint of the divergence operator, 
preserving the divergence is equivalent to solving: 
\begin{equation} \label{ph_eq_bc}
   \int_\Omega \uu_h \cdot s_h(\uu_h, P_n^*\Bgrad_h q_h) 
    + \int_\Omega \Tgrad_{h,p_b} p_h \cdot (P_n^* \Bgrad_h q_h) 
    + \nu \int_\Omega \Tcurl_{h,u_t} \uu_h  \Bcurl_h (P_n^* \Bgrad_h q_h)=0
\end{equation}
for all $q_h \in V_h^2$.

Note that this system is symmetric once one develops the 
definition of $\Tgrad_{h,p_b}$.

We can also see that this discretization of the boundary for the pressure correspond to a homogeneous Neumann boundary condition in $\Gamma_n$ throught the presence of $P_n$ in the definition of $\Tgrad_{h,p_b}$, and a Dirichlet boundary condition on the complement $\Gamma_p$.

\section{Numerical experiments}
\label{sec:Numerics}

We now present numerical results obtained by applying our scheme 
\eqref{eqn:iterCN} to the conforming and broken splines spaces 
described in \cite{buffa2011isogeometric,gucclu2022broken}. 
For the conforming case, the spaces are built via a tensor product of 1D splines of maximal regularity, with anisotropic polynomial degree for the discrete $H(\Div)$ space to ensure the exact sequence property. More precisely, we consider the following discrete de Rham sequence: 
\begin{equation} \label{Discrete_deRham}
  \xymatrix{
   S_{\sf{p}+1} \bigotimes S_{\sf{p}+1} \ar[r]^-{\curl} & 
   (S_{\sf{p}+1} \bigotimes S_{\sf{p}}) \times (S_{\sf{p}} \bigotimes S_{\sf{p}+1}) \ar[r]^-{\Div} & S_{\sf{p}} \bigotimes S_{\sf{p}} .\\
  }
\end{equation}

Where $\sf{p}$ denotes the polynomial degree. We also denote $\sf{n_c}$ the number of cells used to discretize the spaces (cells in which the function are degree $\sf{p} / \sf{p+1}$ polynomial with corresponding $C^{\sf{p-1}}$ or $C^{\sf{p}}$ smoothness condition between the cells). 

For the non-conforming discretization, we cut the domain in each direction in $\sf{n_p}$ sub-domains, called "patches". In each patch we use the previously described spline spaces. The underlying conforming space is recovered by imposing only a $C^0$ constraint in the necessary direction, that is both direction for the discrete $H^1$ space and in the direction normal to the interface of the patches in the discrete $H(\Div)$ space (no conformity condition are needed for the discrete $L^2$ space). Note that the conforming space for the non-conforming discretization is not the same as the space used for the conforming discretization, due to this lower smoothness constraint. The conforming discretization corresponds therefore to a non-conforming discretization with $\sf{n_p}=1$.

We will study the accuracy of our scheme, 
the numerical conservation of the invariants claimed 
in Section~\ref{sec:Time_discretization}
and we will compare the results obtained with the conforming and non-conforming spaces.
Our implementation uses the Psydac library \cite{guclu2022psydac}. 

\subsection{Taylor-Green Vortex}
\begin{figure}
    \centering
    \includegraphics[width=0.49\textwidth]{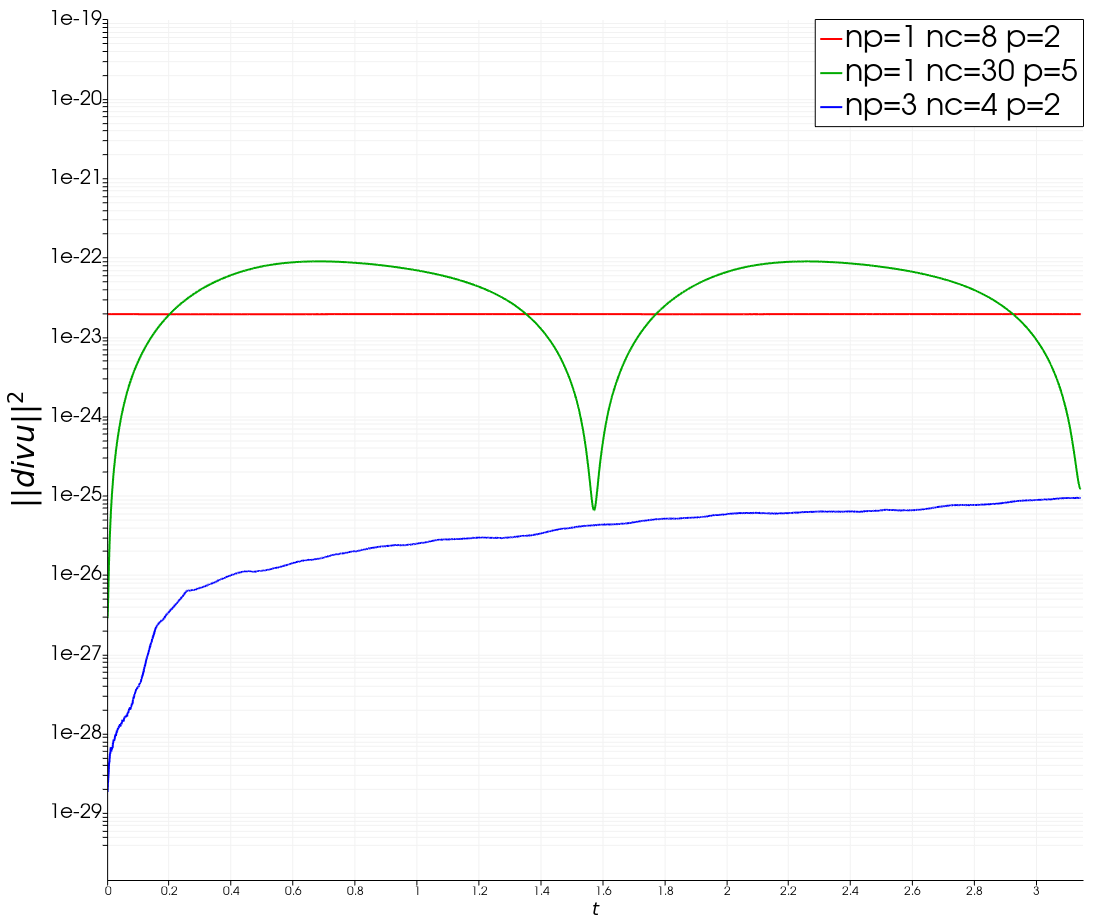}
    \includegraphics[width=0.49\textwidth]{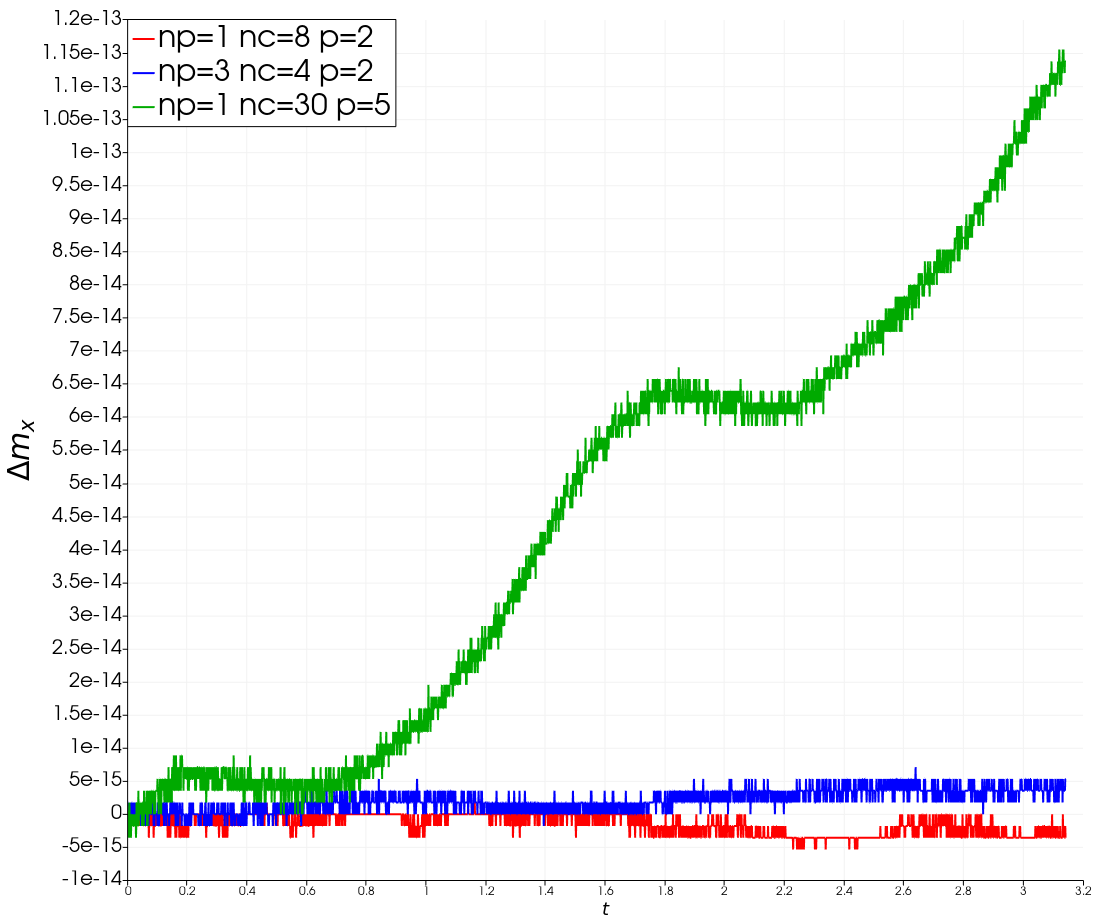}
    \includegraphics[width=0.49\textwidth]{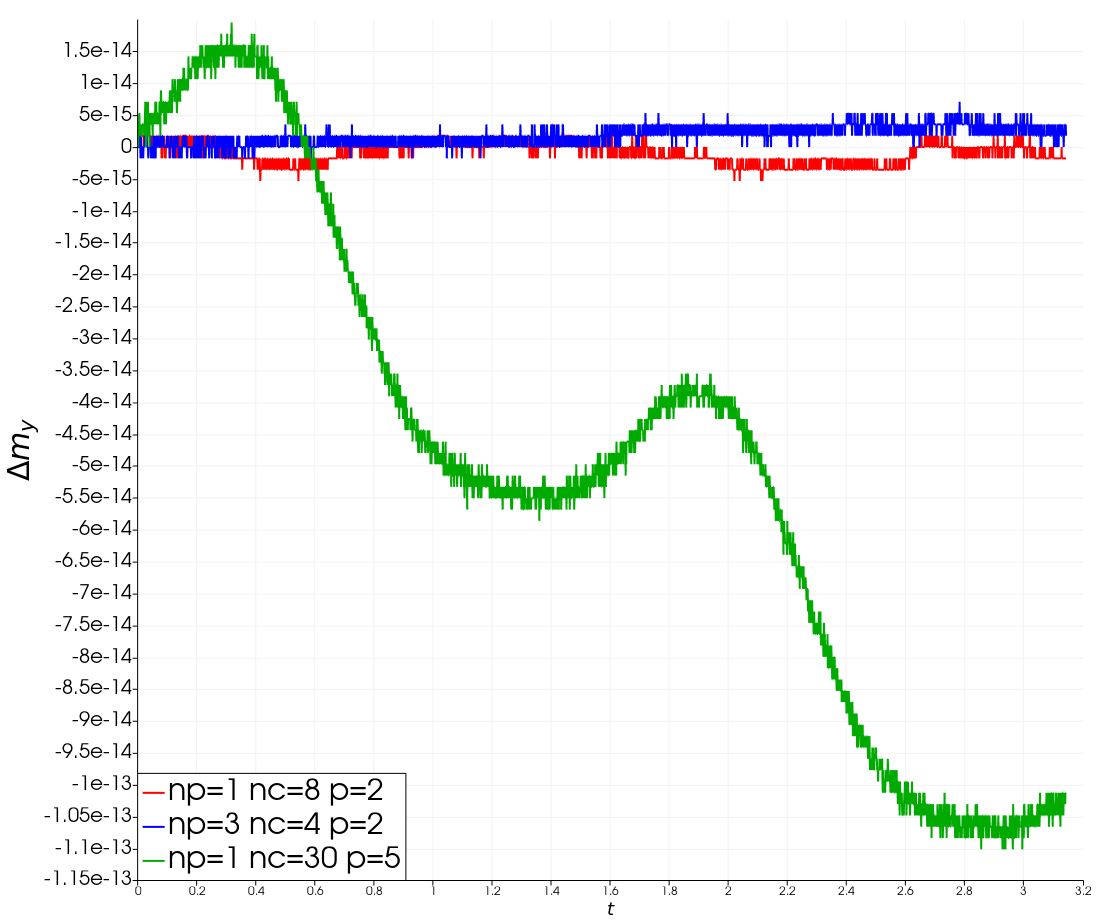}
    \includegraphics[width=0.49\textwidth]{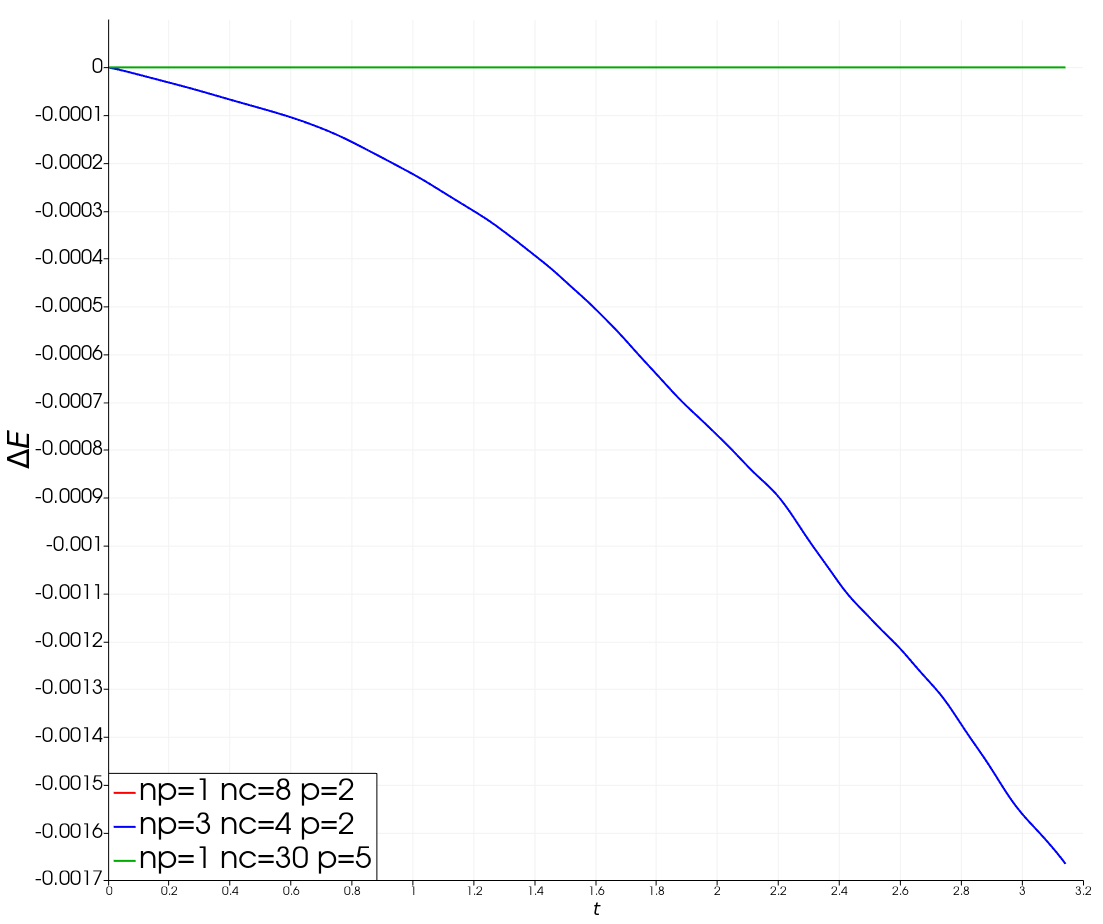}
    \caption{Evolution of the invariants for differents discretisations ($\sf{n_p}$ is the number of patches, $\sf{p}$ the maximum polynomial degree and $\sf{n_c}$ the number of cells per patch}
    \label{fig:invariants}
\end{figure}

\begin{figure}
    \centering
    \begin{subfigure}[b]{0.49\textwidth}
    \centering
    \includegraphics[width=\textwidth]{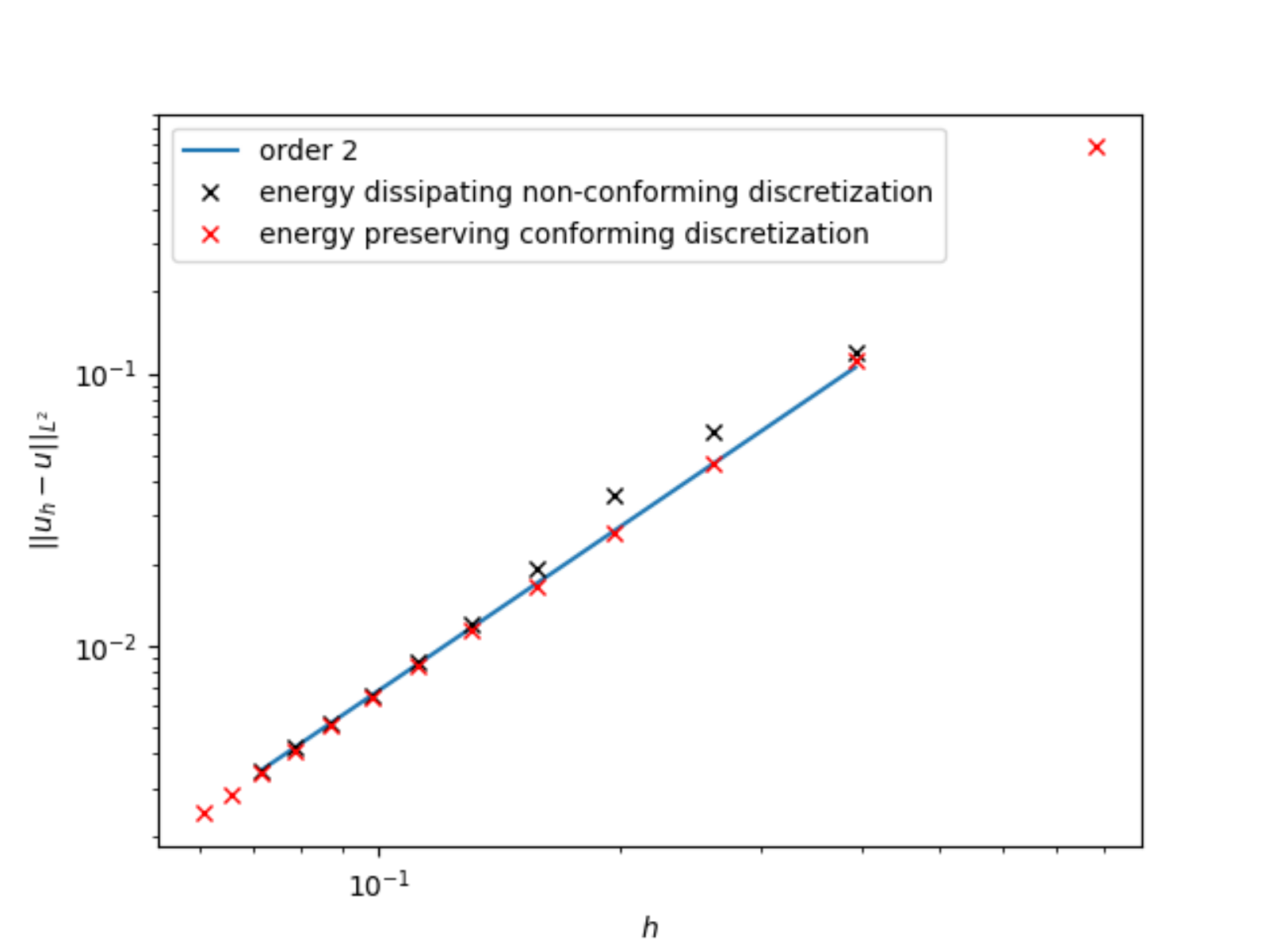}
    \caption{\sf{p}=1}
    \end{subfigure}
    \hfill
    \begin{subfigure}[b]{0.49\textwidth}
    \centering
    \includegraphics[width=\textwidth]{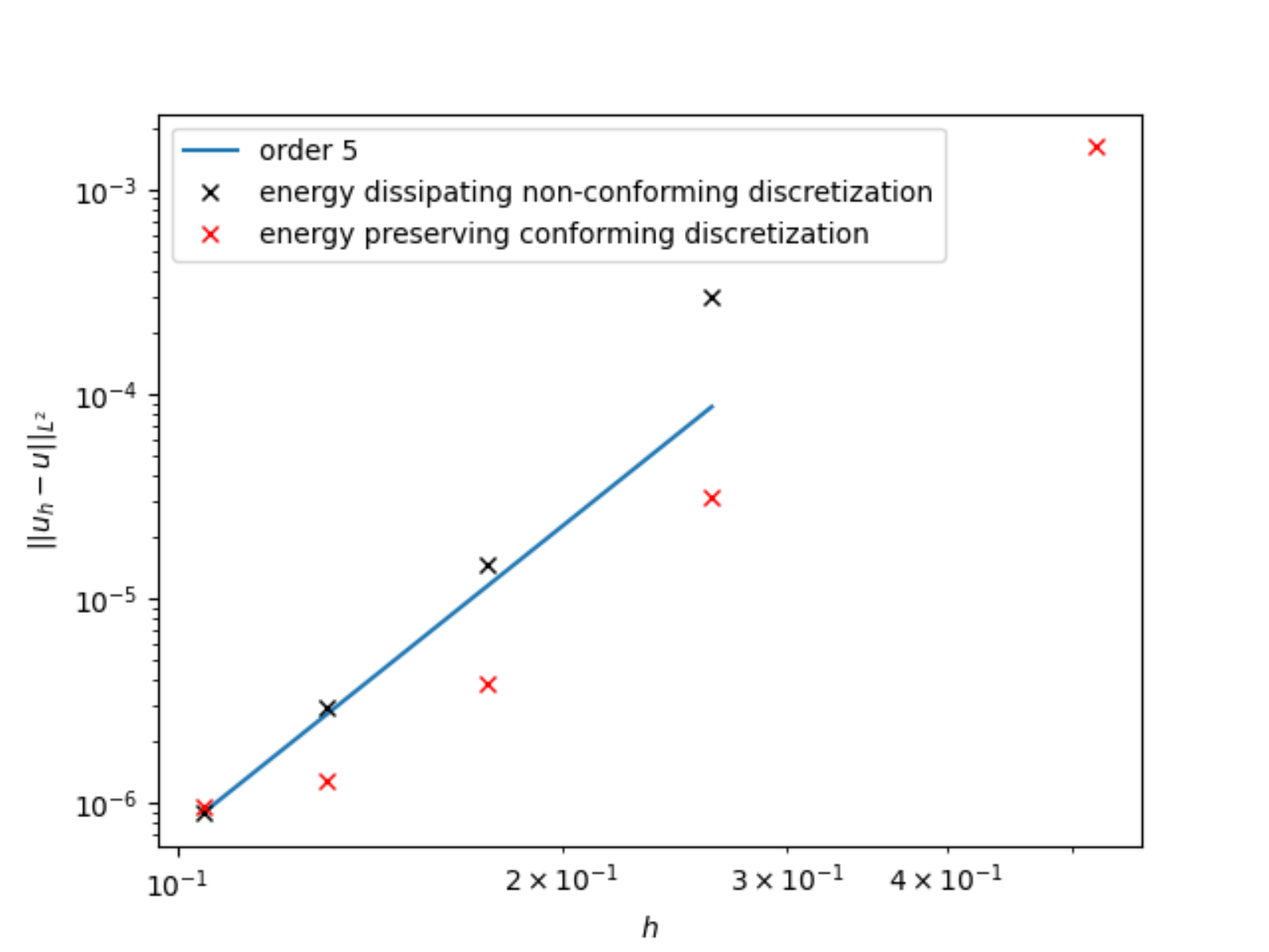}
    \caption{\sf{p}=4}
    \end{subfigure}
    \caption{Evolution of the error for conforming and non-conforming discretizations}
    \label{fig:convergence_Taylor}
\end{figure}

Our first test is a moving Taylor Green vortex on a periodic domain. 
Here the domain is a square torus of period $\pi$ and the solution is given by: 
\begin{equation}
\label{eqn:Taylor-Green}
\uu(x,y;t)=(1-2\cos(2(x-t))\sin(2(y-t)) ,1+2\cos(2(y-t))\sin(2(x-t))) ~,
\end{equation}
with viscosity $\nu=0$. Our first objective is to verify the conservation properties, 
as well as the high order accuracy of our scheme, in a purely advective case. 
The tolerance in the iterative solver is set to $10^{-8}$, and in the non-conforming case 
the dissipation parameter is chosen ad-hoc to be $\alpha = 1000$. The simulation is run until $t_f=1$ with a constant time step of $\Dt = 10^{-4}$, small enough so that 
spatial discretization errors are dominating.
\Cref{fig:invariants} shows the evolution of the invariants for various conforming and 
non-conforming simulations.
One can clearly see that they are well preserved independently of the grid and of the conformity of the simulation (except for the energy which is dissipated in the non-conforming case). \Cref{fig:convergence_Taylor} shows the $L^2$ errors as the meshes are refined, for different polynomial orders: in both the conforming and non-conforming cases we observe a high order convergence corresponding to
the optimal rate.

\subsection{Poiseuille flow}
Our second test case is the Poiseuille flow in a channel. Here the domain is the square $[0,\pi]^2$ and we impose the non-slip boundary conditions on the left and right boundaries, that is $\uu=0$ on $\{0\} \times [0,\pi] \cup \{\pi\} \times [0,\pi]$. 
We impose the pressure at the top and at the bottom as well as a normal velocity constraint, that is $\uu \times \nn =0 $ on $[0,\pi] \times \{0\} \cup [0,\pi] \times \{\pi\}$ 
and $p=\pi/2$ on $[0,\pi] \times \{\pi\}$ and $p=-\pi/2$ on $[0,\pi] \times \{0\}$. 
With these settings there is a steady solution, given by 

\begin{equation}
\label{eqn:Poiseuille}
\uu(x,y)=(0,\frac{\pi}{\nu} \frac{x (x-\pi)}{2})
\end{equation}

\begin{figure}
    \centering
    \includegraphics[width=13cm]{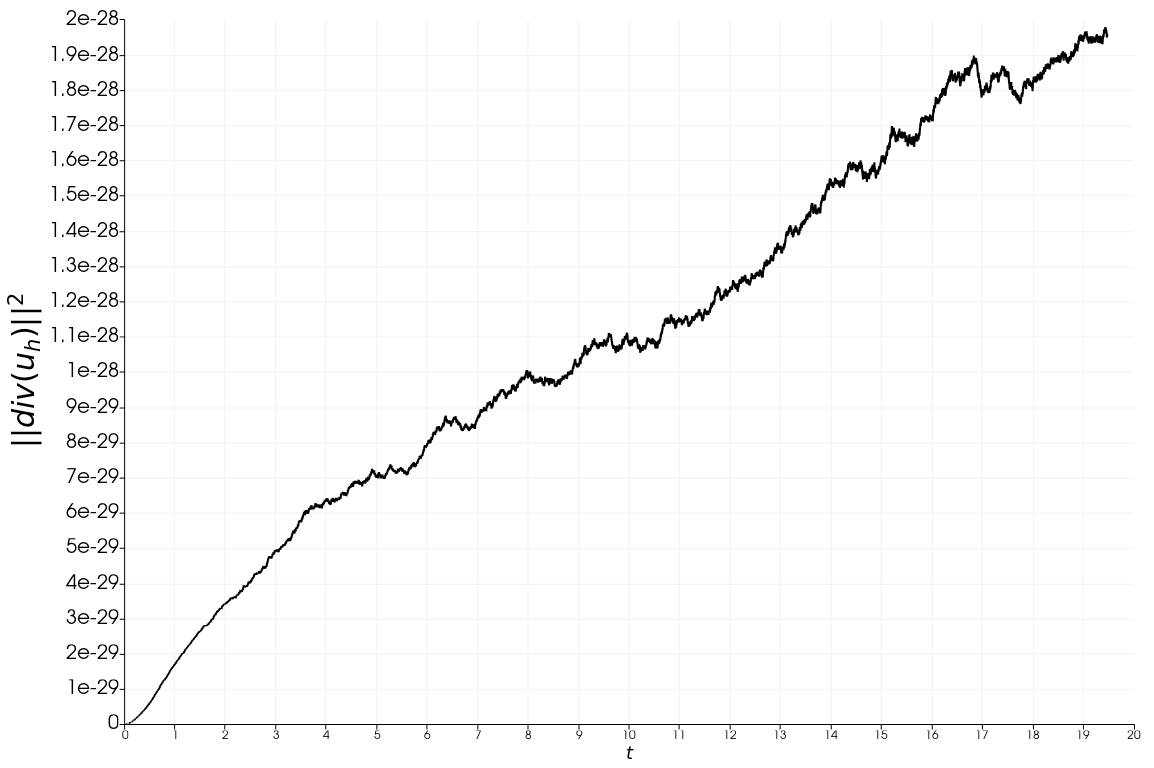}
    \caption{Evolution of the divergence for a simulation of the Poiseuille flow}
    \label{fig:evol_div_Poiseuille}
\end{figure}

In this case there is no more conservation of the energy and momentum (due to the presence of boundaries), however we still have preservation of the divergence as shown in \cref{fig:evol_div_Poiseuille}, where we can see that our handling of the boundary conditions does not affect the preservation of mass, as claimed in \cref{sec:pres_eq_bnd}. For $\sf{p}$ greater or equal to 2, the exact solution \cref{eqn:Poiseuille} is in the discretization space. As expected the error is then comparable to the error done in the non-linear system, for example with a discretization of $3 \times 3$ patches, $4\times 4$ cells per patch and the linear solver tolerance set to $10^{-8}$, the final error is $3.37 \times 10^{-7}$ proving that our scheme is able to approximate the solution optimally.

\subsection{Lid-driven cavity}
We now move on to the lid-driven cavity problem \cite{GGS82}  which is widely used to test the consistency and accuracy of the implementation of an incompressible or low Mach flow solver, see e.g. \cite{Bassi2006,Boscheri2021,Peshkov2021,Fambri2023}. Indeed, this test provides a 2d example where all the terms of the PDE are equally and strongly important but, in contrast, wall-boundary conditions are specified on all boundaries. 

In particular, at the upper wall a positive tangential velocity is imposed through the numerical discretization of the viscous terms. This makes the problem difficult 
since the prescribed boundary conditions are \emph{discontinuous} at the upper corners of the cavity, with opposite pressure peaks. 

Specifically, here the domain is the square $\Omega= [0,1]^2$ and we impose 
a no-slip boundary condition on the left, right and bottom boundaries: 
$\uu=(0,0)$ in $\{0\} \times [0,1] \cup \{1\} \times [0,1] \cup [0,1] \times \{0\}$. 
On the top boundary we impose a viscous moving plate condition: 
$\uu = (1,0)$ on $[0,1] \times \{1\}$. 
The fluid viscosity is $\nu=10^{-2}$, which correspond to a Reynolds number of $Re=100$. For this test, the discretization is done by choosing $10\times 10$ patches, splines with polynomial degree $\sf{p}=2$ and $4\times 4$ cells per patch. This numerical set-up corresponds to $7\times 6$ degrees of freedom per patch for the horizontal velocity $u_x$, $6\times 7$ for the vertical velocity $u_y$ and $6\times 6$ for the pressure $p$.

\begin{figure}
    \centering
    \includegraphics[width=0.7\textwidth,trim = 6cm 0 6cm 0,clip]{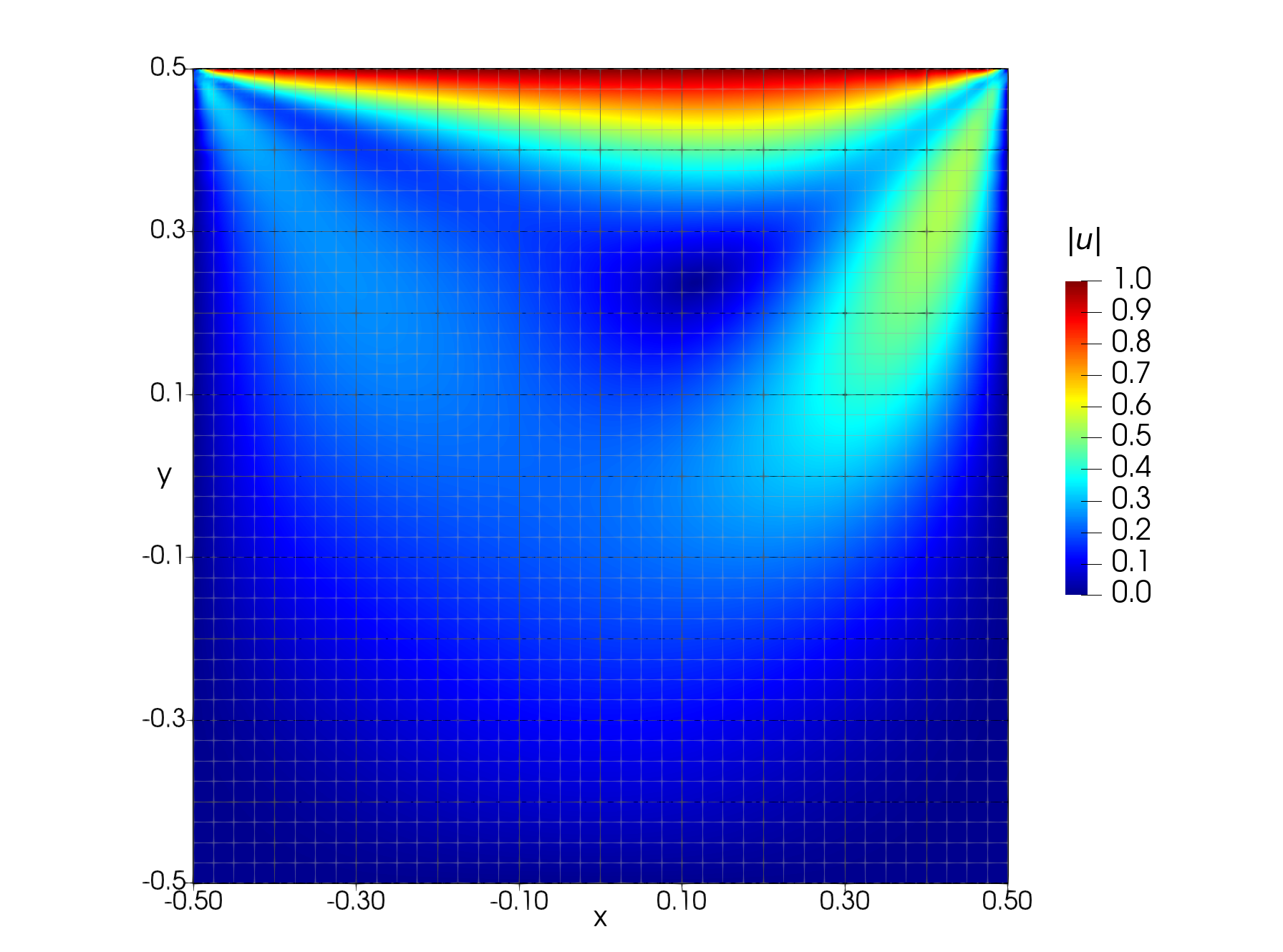}
    \caption{Norm of the velocity of the solution to the lid driven cavity problem for $\nu=10^{-2}$, with $10\times 10$ patches, degree $\sf{p}=2$ and $4\times 4$ cells per patch. This numerical set-up corresponds to $7 \times 6$ degrees of freedom per patch.}
    \label{fig:sol_LDC}
\end{figure}

\begin{figure}
    \centering
    \includegraphics[width=0.7\textwidth]{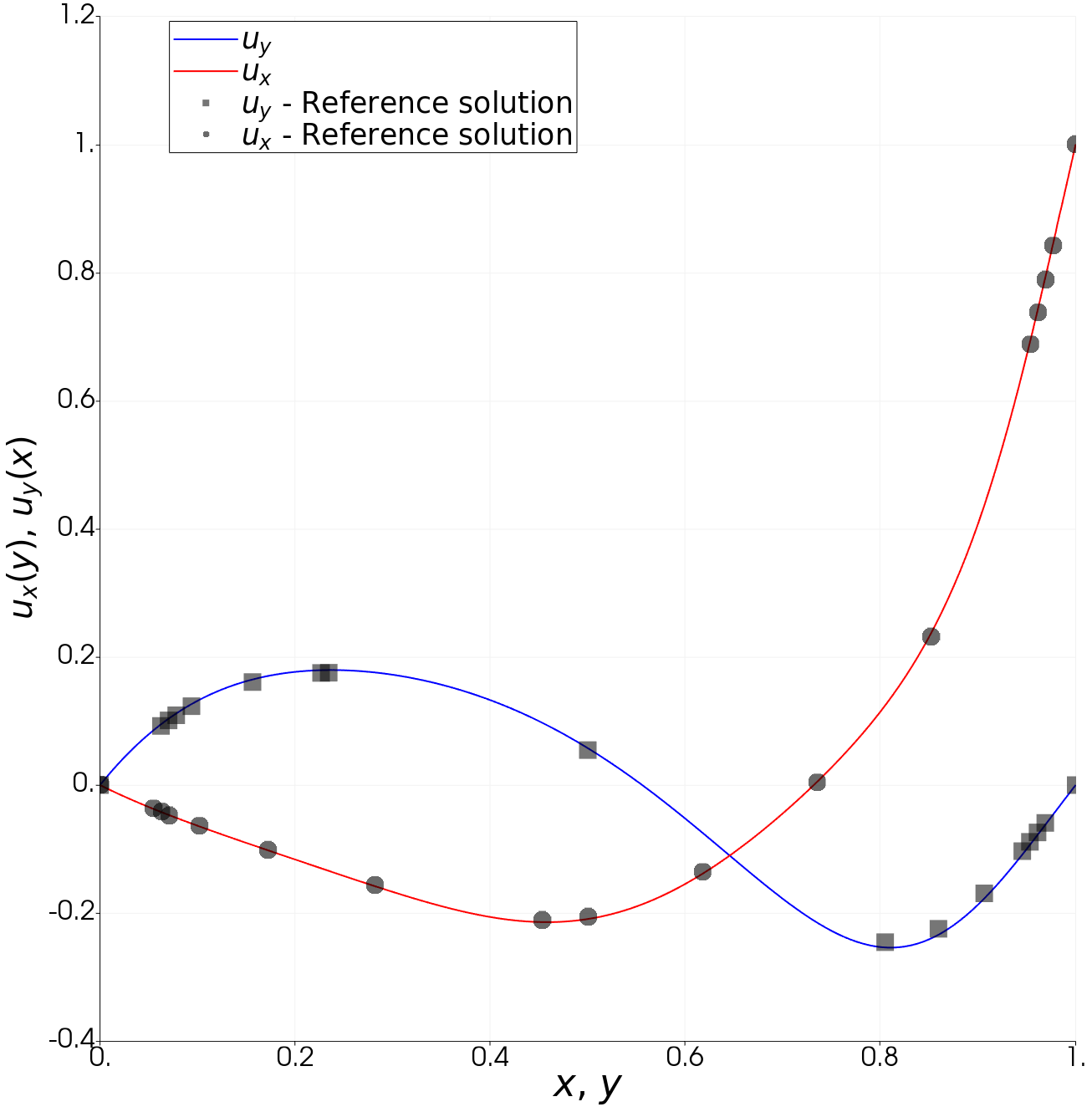}
    \caption{Cuts along vertical ($x=0$) and horizontal ($y=0$) axes for the solution $u_x(0,y)$ and $u_y(x,0)$, respectively, for the lid-driven cavity problem. Reference
    values are taken from Ghia et al. \cite{GGS82}.}
    \label{fig:compar_LDC}
\end{figure}

\Cref{fig:sol_LDC} shows the computed solution to this problem and \cref{fig:compar_LDC} 
compares our solution and a reference one taken from \cite{GGS82}, along vertical ($x=0$)
and horizontal ($y=0$) cuts. We can see that our new method is able to reproduce well this flow, proving the robustness of our discrete Poisson solver and the correctness of the discrete boundary conditions.

\subsection{Blasius boundary layer}
\begin{figure}
    \centering
    \includegraphics[width=\textwidth,trim = 0 8cm 0 8cm,clip]{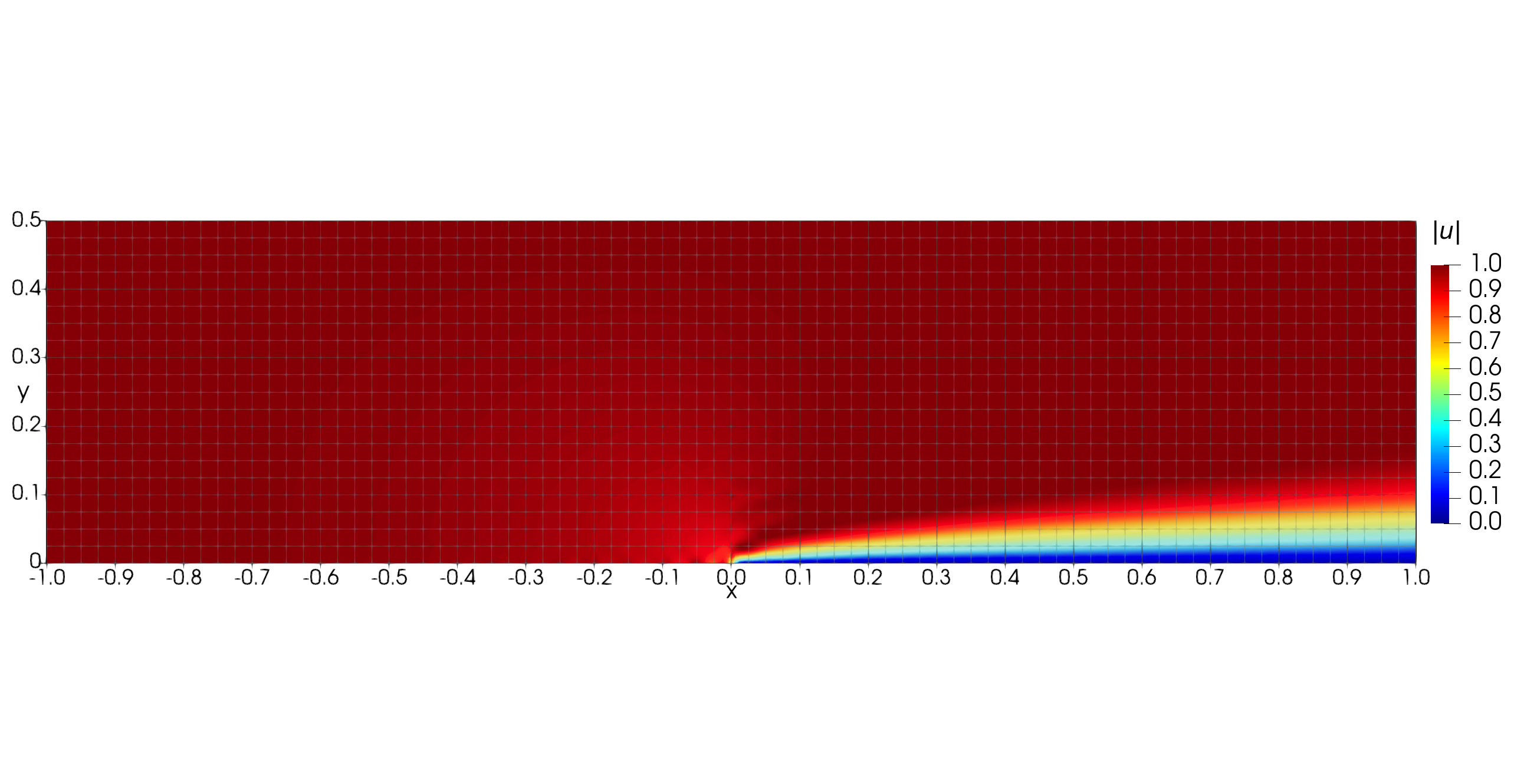}
    \caption{Norm of the velocity for the solution to the Blasius problem for $\nu=10^{-3}$, with $20\times 5$ patches, degree $\sf{p}=2$ and $4\times 4$ cells per patch.}
    \label{fig:sol_Blasius}
\end{figure}

\begin{figure}
    \centering
    \includegraphics[width=13cm]{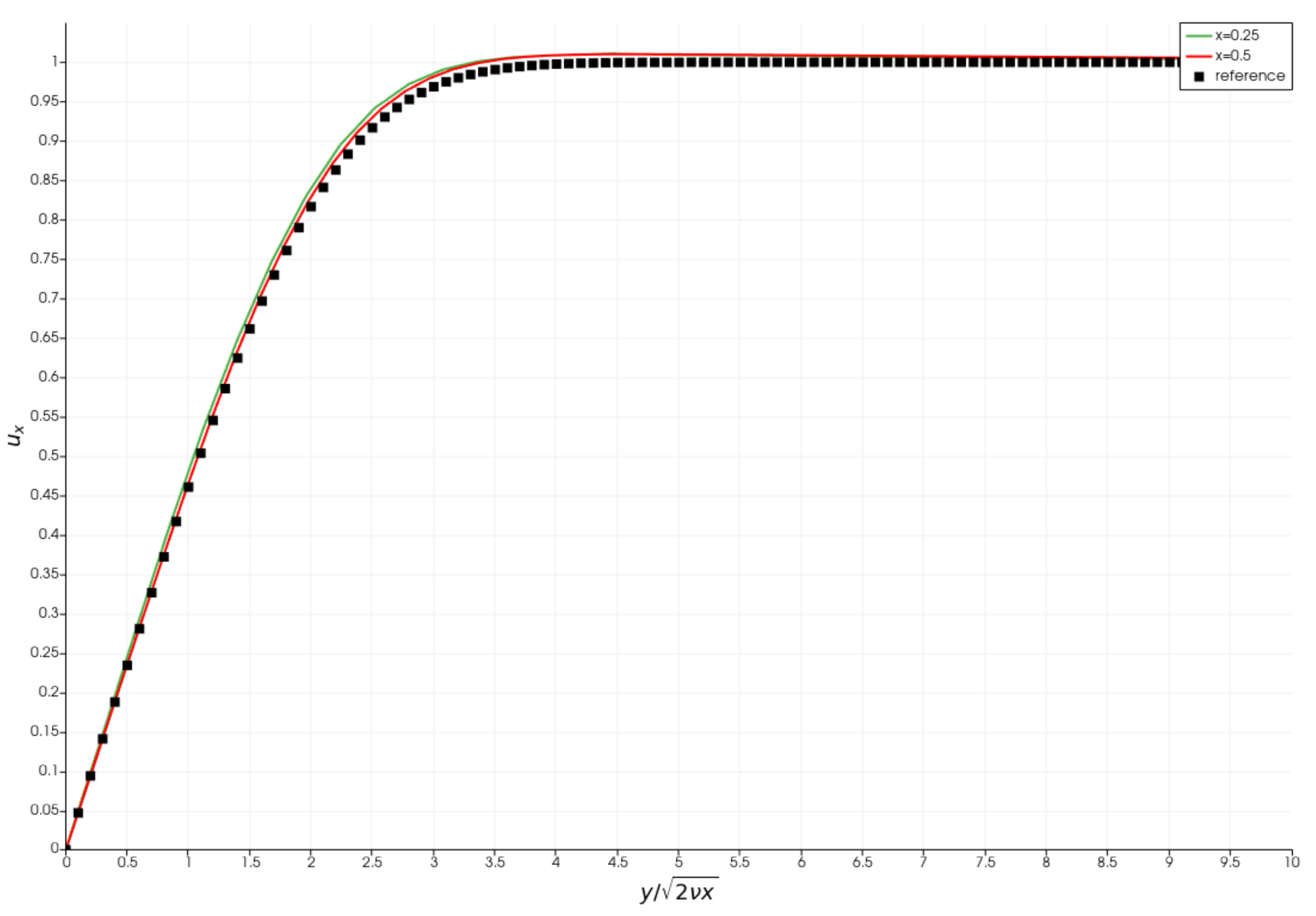}
    \caption{1d cuts along vertical axis for the solution of the Blasius (laminar) boundary layer and comparison with reference.}
    \label{fig:compar_Blasius}
\end{figure}

Our last test is the Blasius (laminar) boundary layer, where we investigate the accuracy of our method near the boundaries in the treatment of the viscosity. The physical domain is a rectangle $\Omega=[-1,1]\times [0,0.5]$, and the viscosity is $\nu=10^{-3}$.
We impose for boundary conditions a constant inflow on the left, that is $\uu=(1,0)$ on $\{-1\}\times [0,0.5]$, a slip boundary condition on the left part of the bottom ($\uu \cdot \nn = 0$ on $[-1,0]\times \{0\}$) and no-slip boundary condition on the right part ($\uu = (0,0)$ on $[0,1]\times \{0\}$). On the top and right boundaries we impose $p = 0$.
Again, similarly to the lid-driven cavity problem, strong pressure gradients are generated where the flow first meets the wall at $(x,y)=(0,0)$, 
and this make any high-order finite-element scheme exposed to oscillations.  
The numerical solution is obtained after choosing   $20\times 5$ patches, local splines of degree $\sf{p}=2$ and $4\times 4$ cells per patch. This numerical set-up corresponds to $7\times6$ degrees of freedom per patch for the horizontal velocity $u_x$, $6\times7$ for the vertical velocity $u_y$ and $6\times6$ for the pressure $p$.  The dissipation parameter used for
conformity errors is $\alpha = 100$, and the solver tolerance is set to $10^{-8}$.
We can see on \cref{fig:sol_Blasius} the formation of the boundary layer. \Cref {fig:compar_Blasius} compares two vertical cuts at $x=0.25$ and $x=0.5$  with a 
reference solution taken from \cite{schlichting_boundary-layer_2017}. 
Note that for these cuts, the boundary layer is almost completely solved within a single patch. Despite the coarse spatial-resolution, we can see that our method is able to have accurate solutions near boundaries.

\subsection{Double shear layer} \label{sec:DSL}
\begin{figure}
	\centering
	\includegraphics[width=0.43\linewidth,clip,trim=2cm 2cm 8cm 2cm]{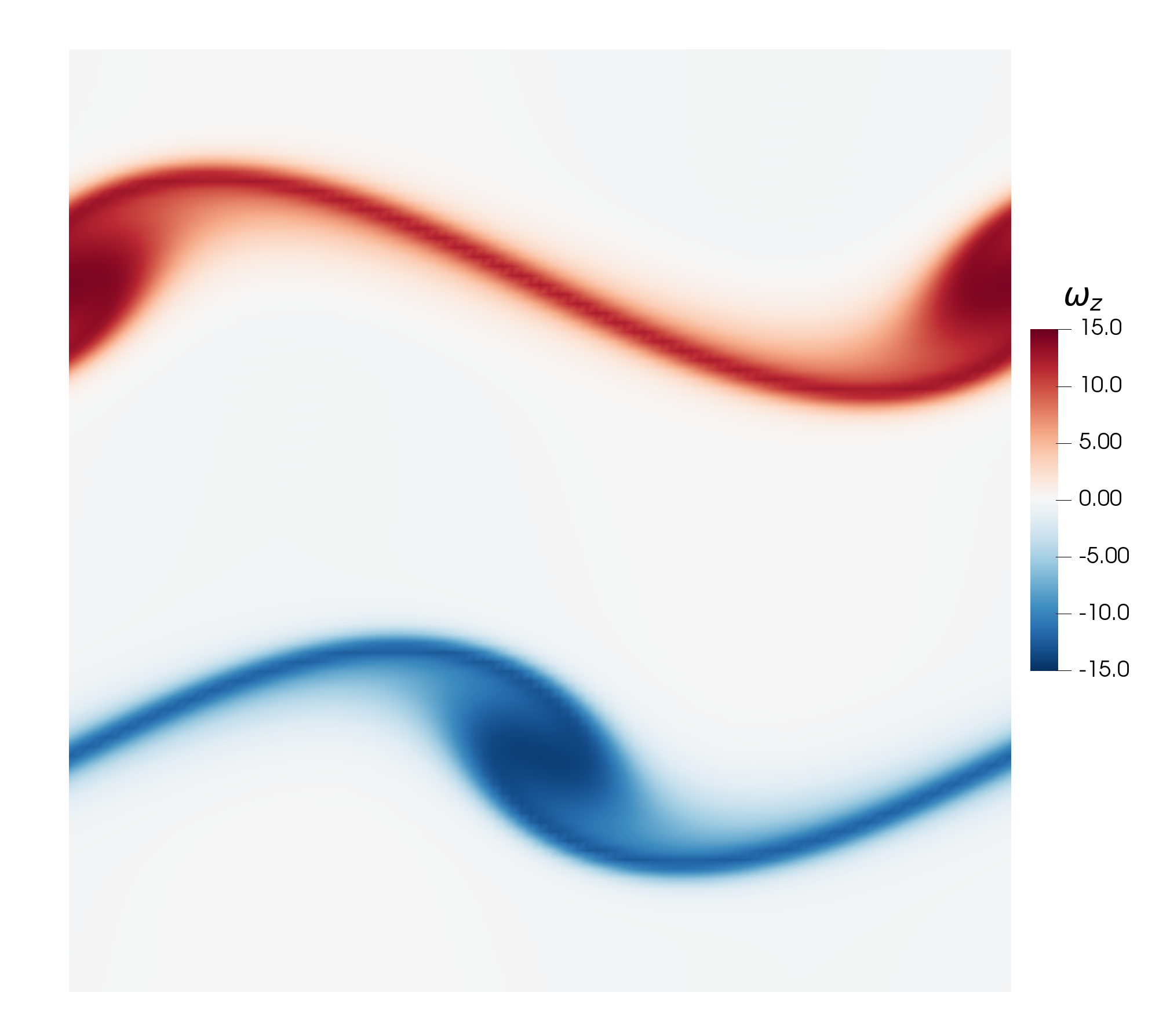} \hspace{0.01\linewidth}	\includegraphics[width=0.43\linewidth,clip,trim=2cm 2cm 8cm 2cm]{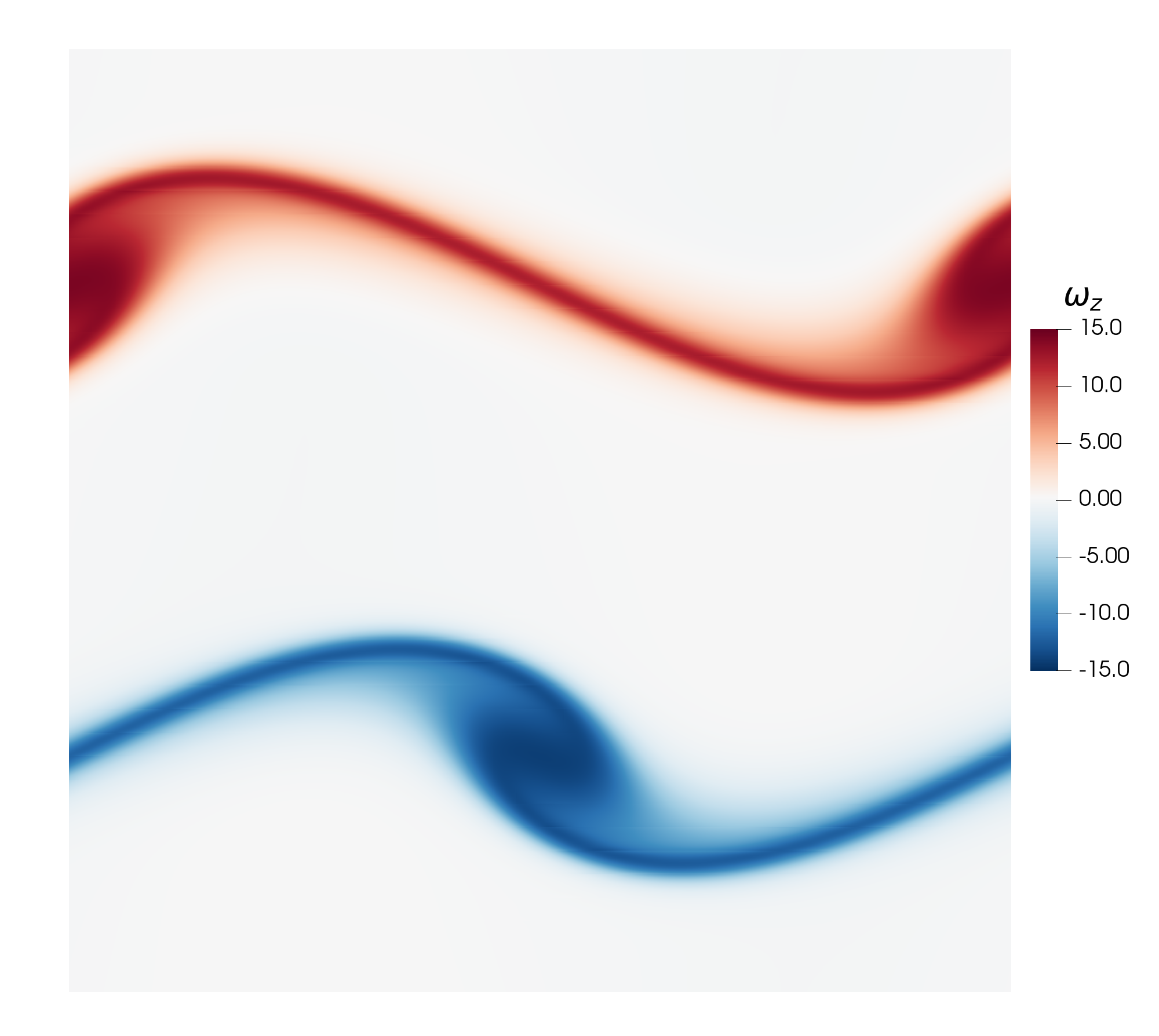}\\
	\includegraphics[width=0.43\linewidth,clip,trim=2cm 2cm 8cm 2cm]{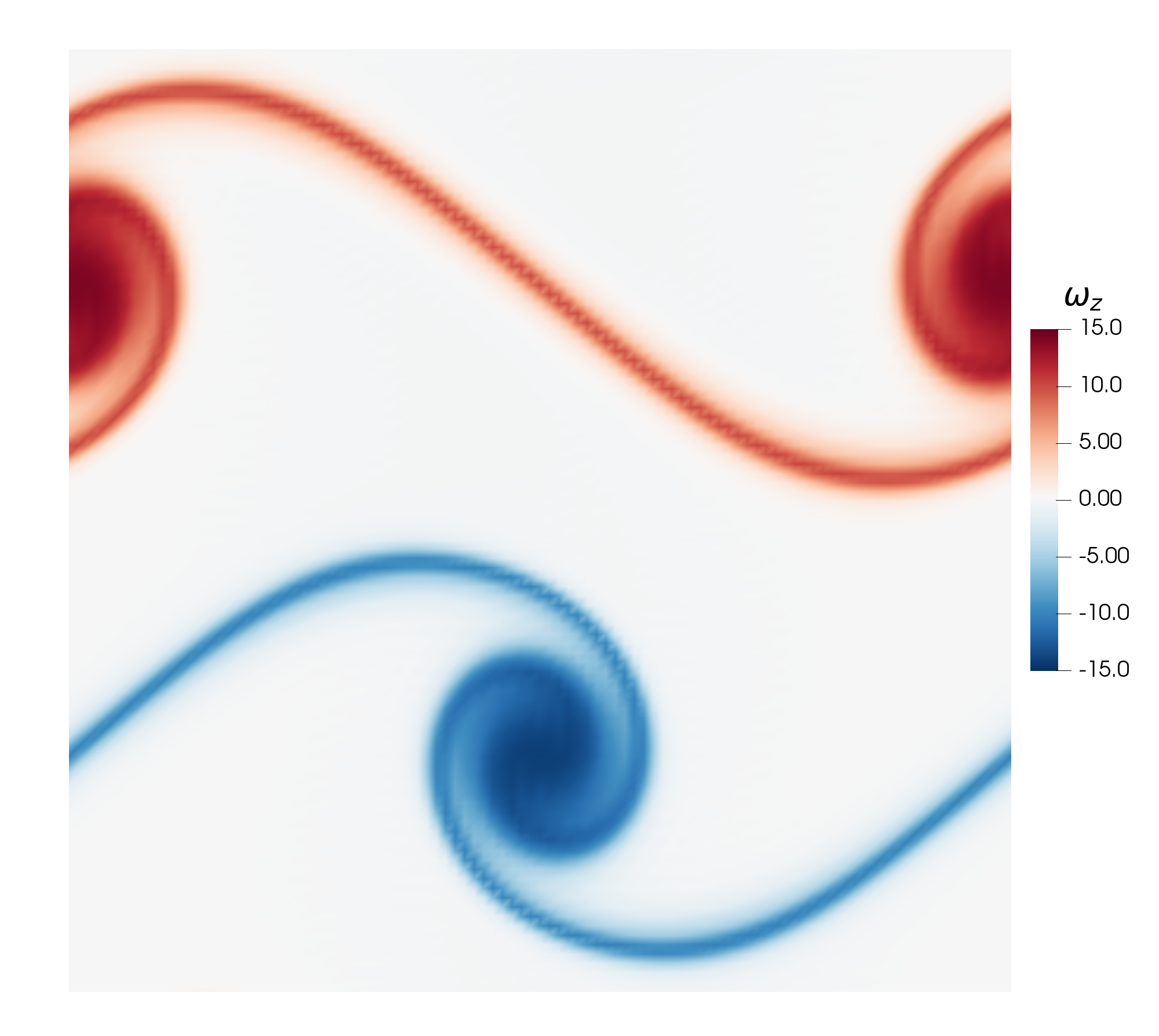}\hspace{0.01\linewidth} 	\includegraphics[width=0.43\linewidth,clip,trim=2cm 2cm 8cm 2cm]{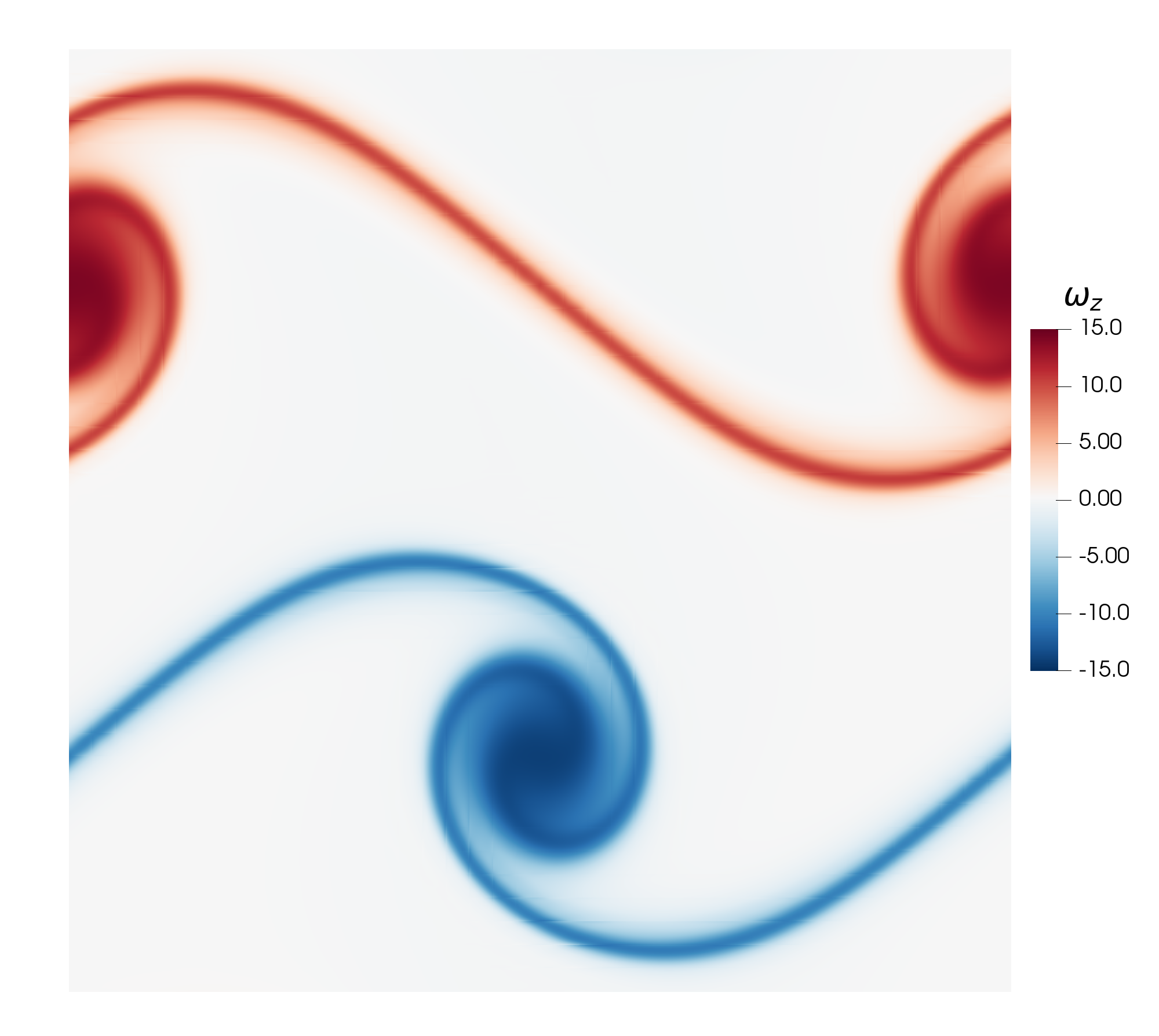}\\
	\includegraphics[width=0.43\linewidth,clip,trim=2cm 2cm 8cm 2cm]{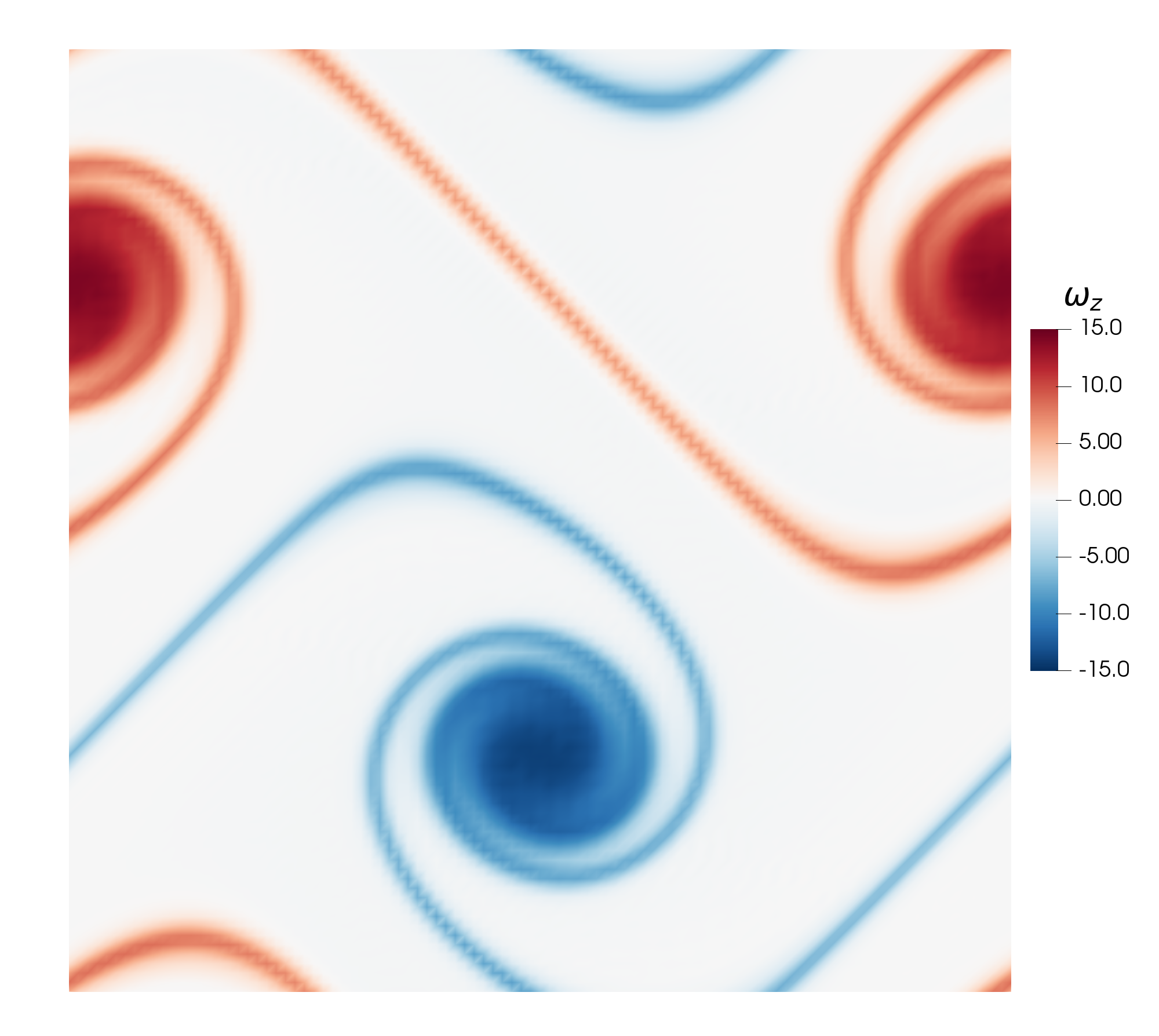}\hspace{0.01\linewidth} 	\includegraphics[width=0.43\linewidth,clip,trim=2cm 2cm 8cm 2cm]{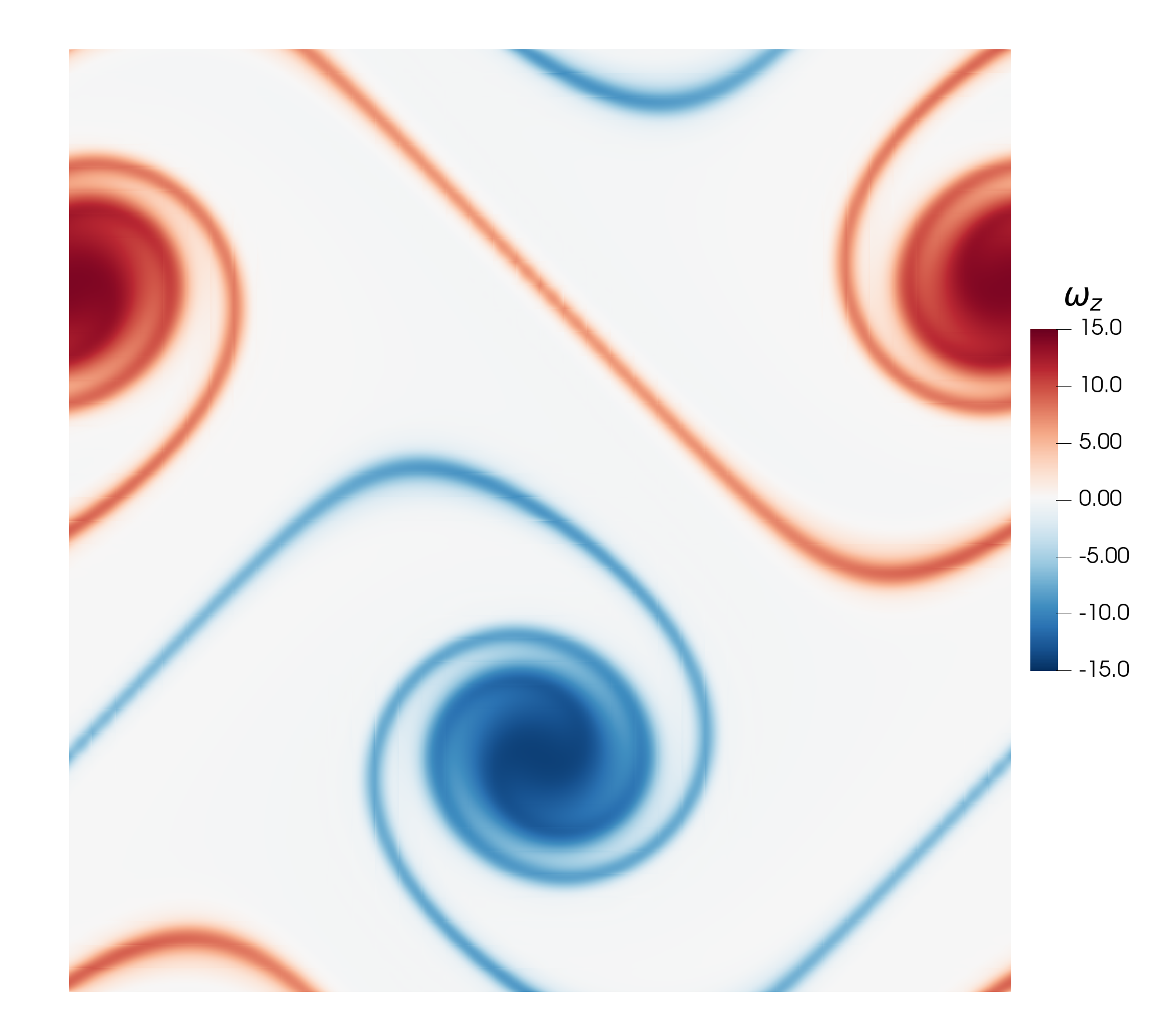}\\
	\includegraphics[width=0.8\linewidth]{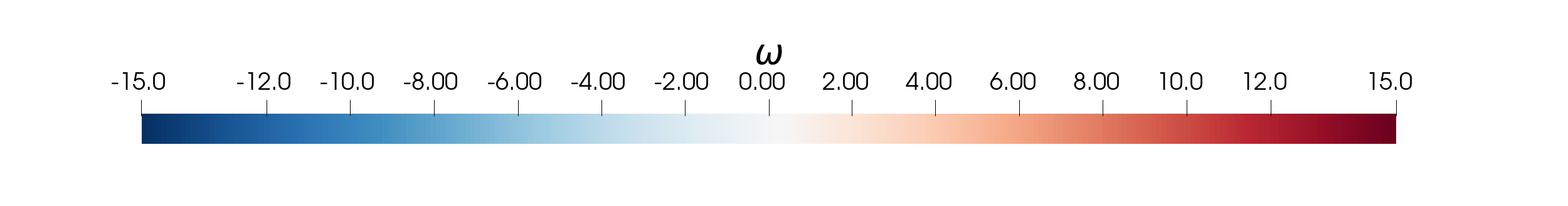}
	\caption{Vorticity for the solution to the double-shear layer problem at times $t=1.5$, $2.0$, $2.5$, from the top to the bottom, for $\nu=0.0002$, with one single periodic patch, degree $\sf{p}=2$ and $100\times 100$ cells. The numerical solution obtained with our conforming FEEC discretization is plotted (left) next to the reference solution, computed at the costs of a 5th order accurate spectral discontinuous Galerkin scheme on a staggered $40\times40$ grid, see \cite{FambriDumbser}.}   \label{fig:DSL1}
\end{figure}

\begin{figure}
	\centering
	\includegraphics[width=0.43\linewidth,clip,trim=2cm 2cm 8cm 2cm]{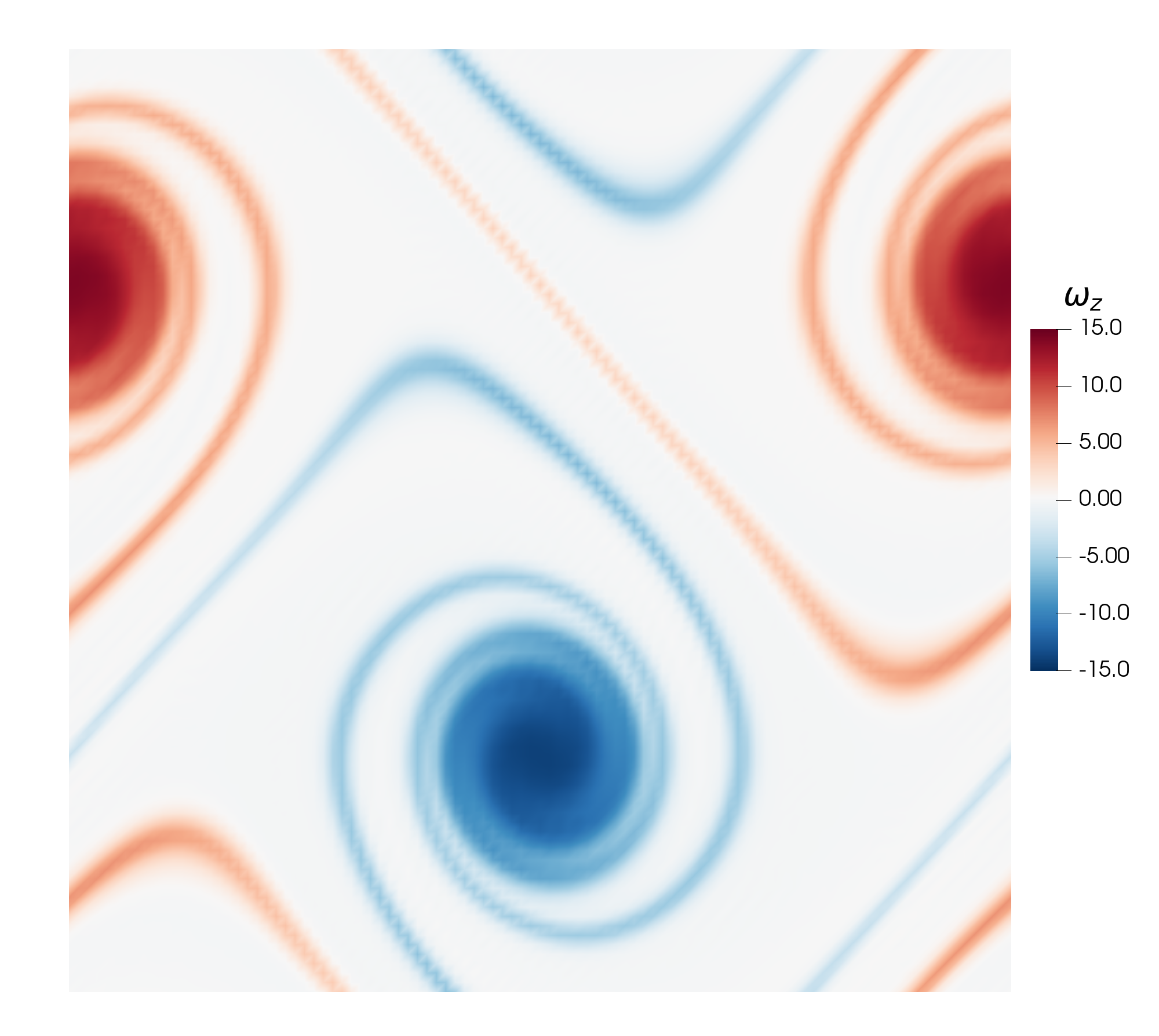} \hspace{0.01\linewidth}	\includegraphics[width=0.43\linewidth,clip,trim=2cm 2cm 8cm 2cm]{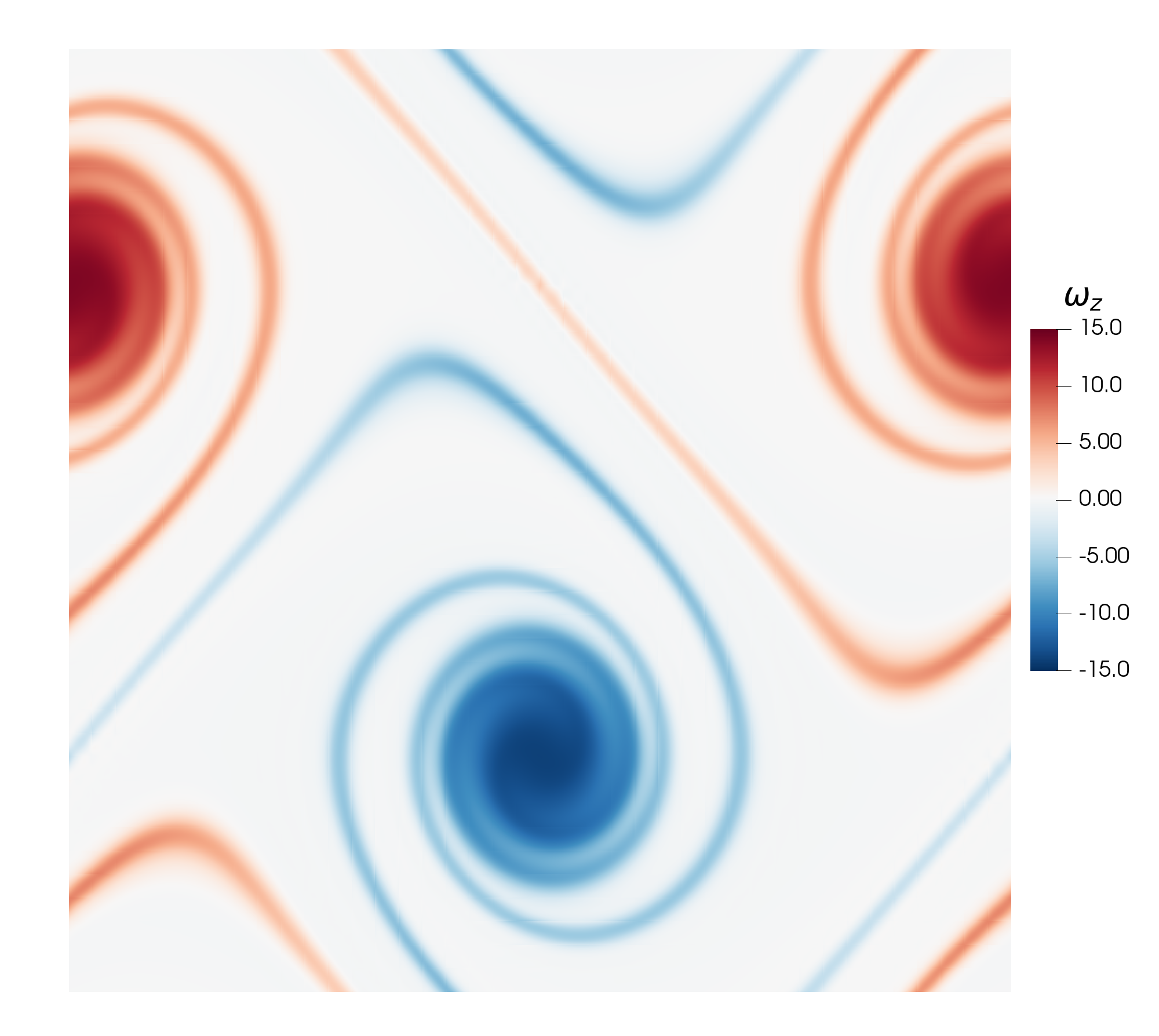}\\
	\includegraphics[width=0.43\linewidth,clip,trim=2cm 2cm 8cm 2cm]{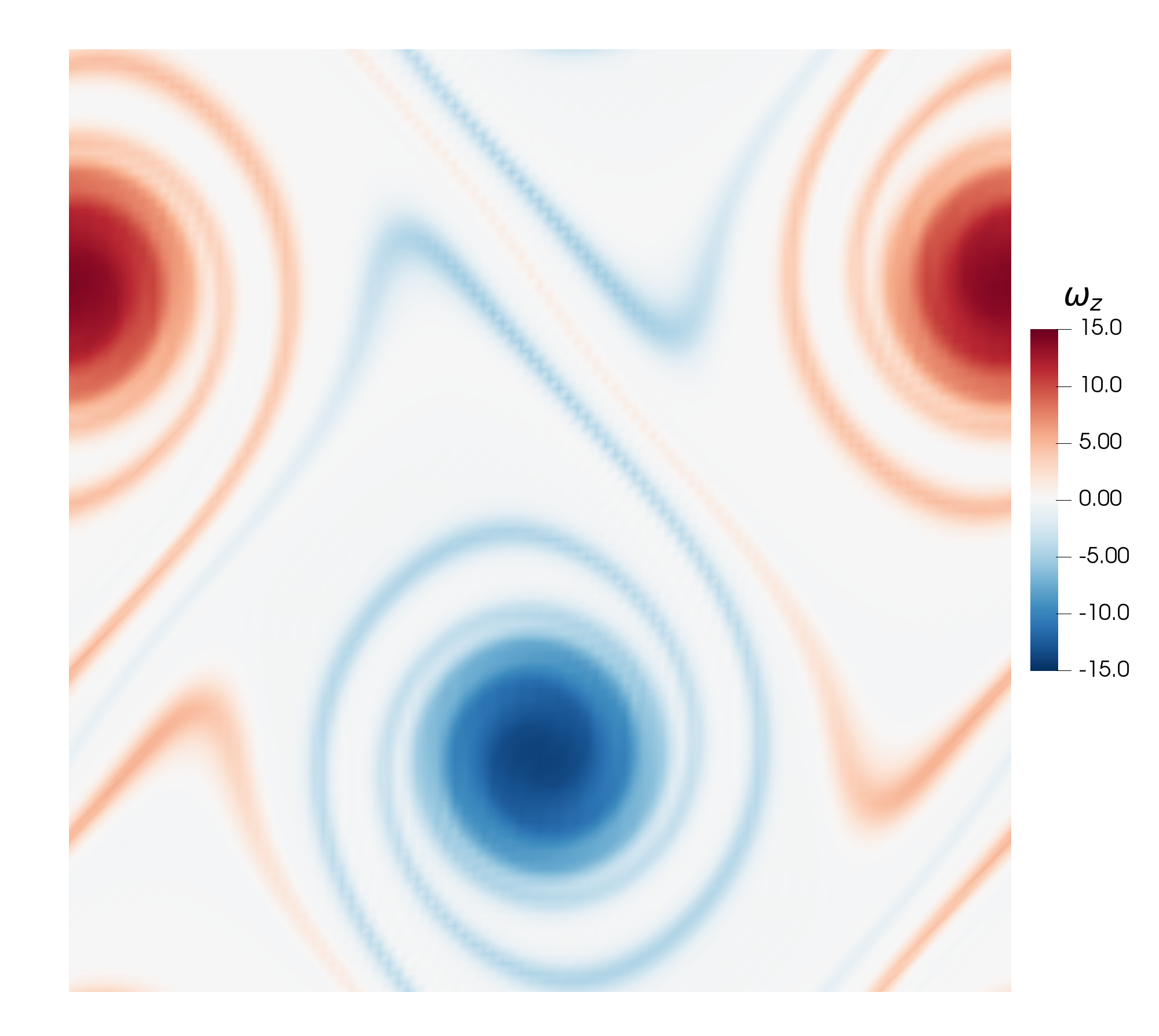}\hspace{0.01\linewidth} 	\includegraphics[width=0.43\linewidth,clip,trim=2cm 2cm 8cm 2cm]{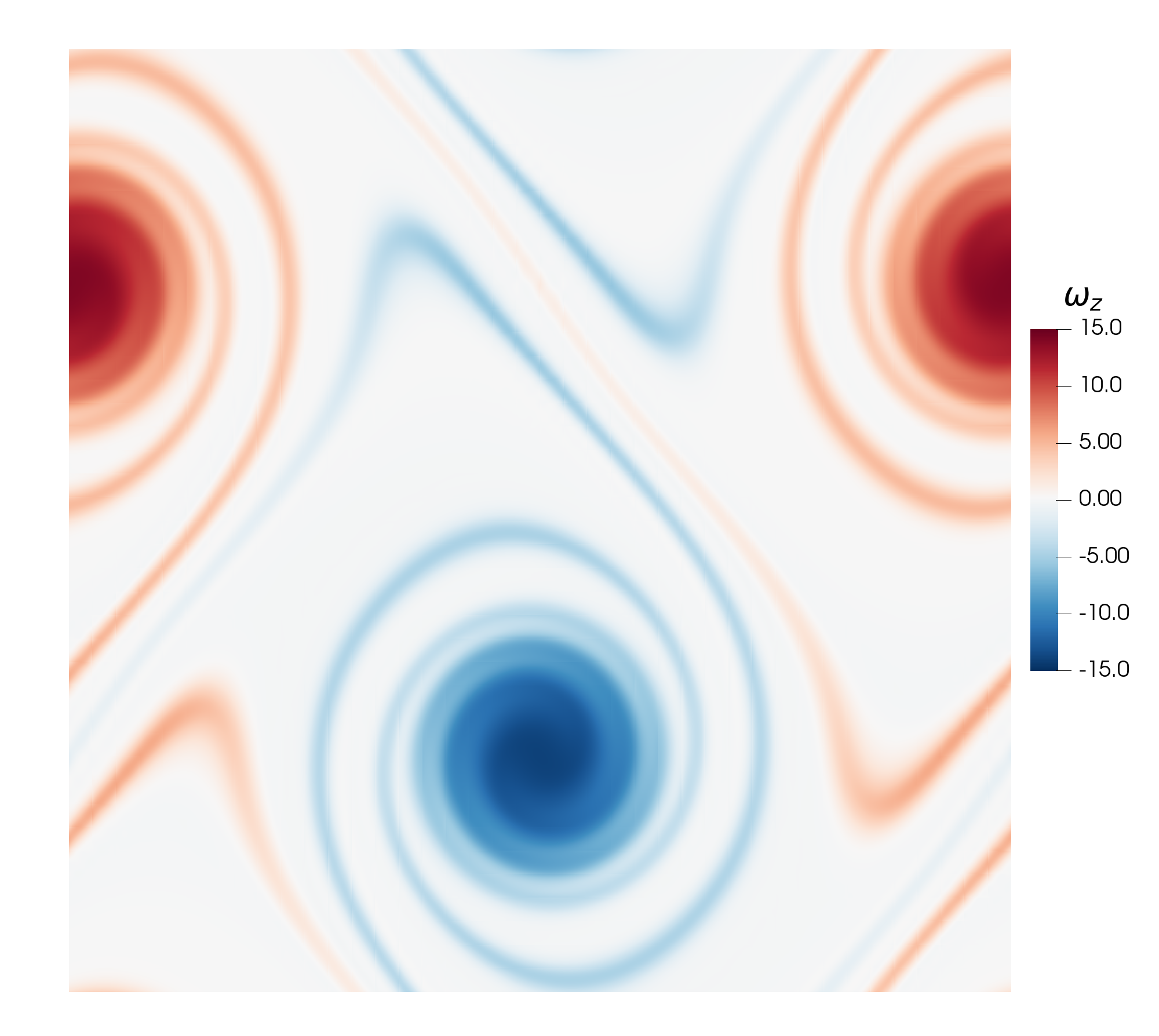}\\
	\includegraphics[width=0.43\linewidth,clip,trim=2cm 2cm 8cm 2cm]{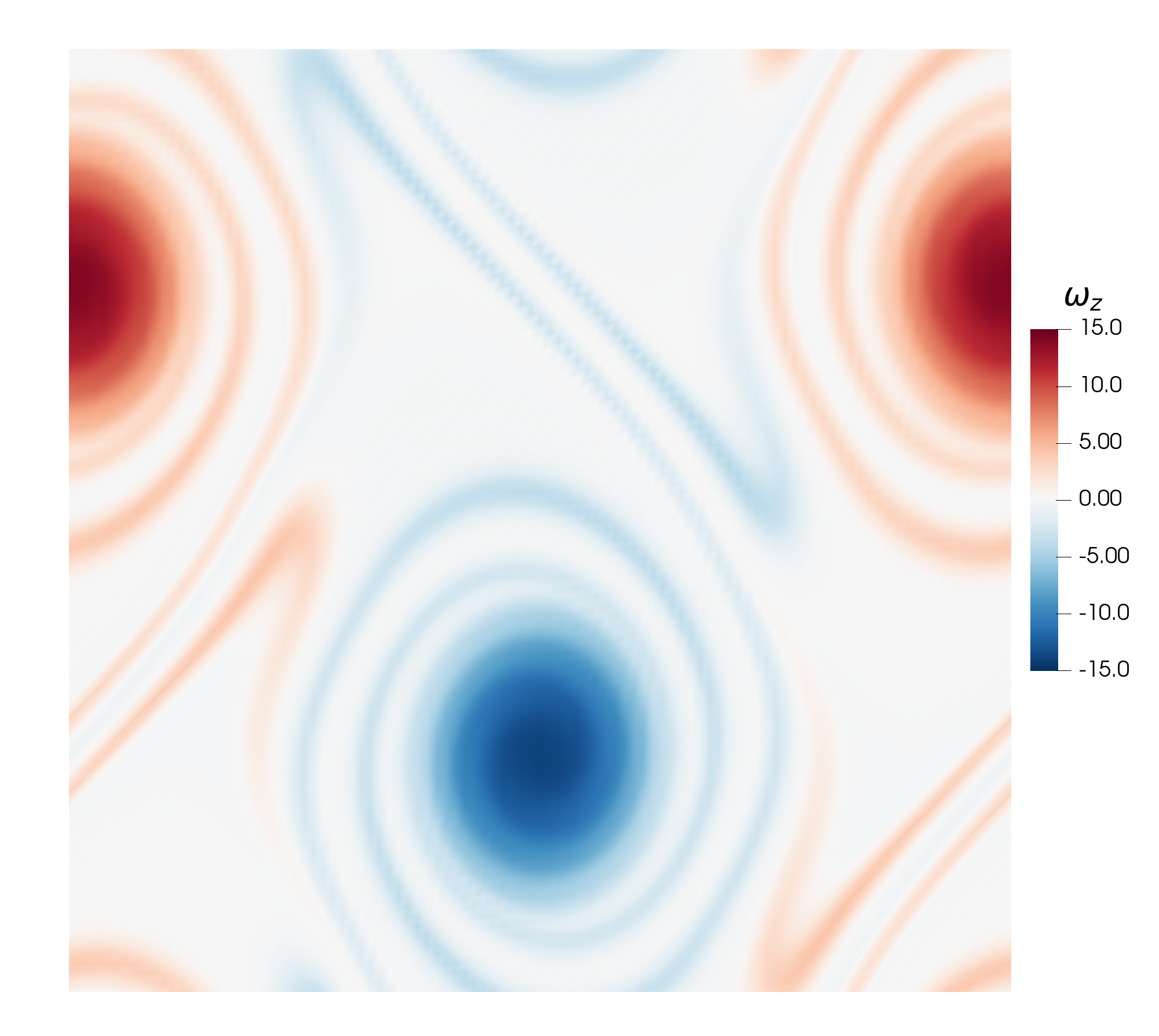}\hspace{0.01\linewidth} 	\includegraphics[width=0.43\linewidth,clip,trim=2cm 2cm 8cm 2cm]{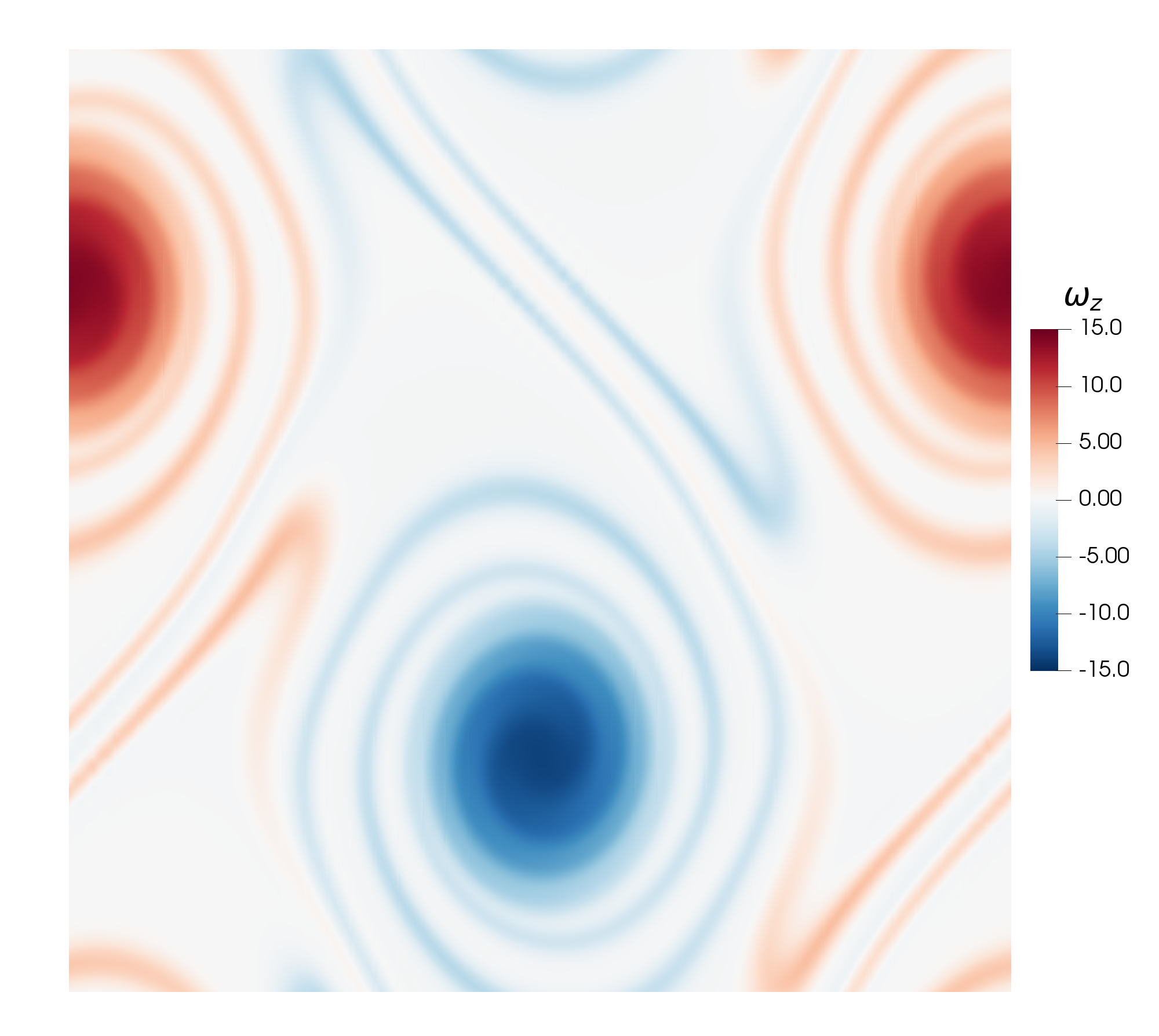}\\
	\includegraphics[width=0.8\linewidth]{DSL_colorlegend}
	\caption{Vorticity for the solution to the double-shear layer problem at times $t=3.0$, $3.5$, $4.0$, from the top to the bottom, for $\nu=0.0002$, with one single periodic patch, degree $\sf{p}=2$ and $100\times 100$ cells. The numerical solution obtained with our conforming FEEC discretization is plotted (left) next to the reference solution, computed at the costs of a 5th order accurate spectral discontinuous Galerkin scheme on a staggered $40\times40$ grid, see \cite{FambriDumbser}.}   \label{fig:DSL2}
	\label{fig.DSL_time}
\end{figure}

Finally we present the numerical results against the so-called double shear layer test \cite{BCG89}. Here, the numerical domain is a periodic box  $\Omega=[-1,1]^2$, and the viscosity is $\nu=0.0002$.
The initial condition consists of
\begin{equation*} 
p=0, \qquad
	\uu(x,y)  = \left(\tanh\Big(\frac{y+0.5}{\delta}\Big)-\tanh\Big(\frac{y-0.5}{\delta}\Big)-1,
		\varepsilon \sin (2\pi x)\right)
\end{equation*}
where $\varepsilon=0.05$ and $\delta = 1/15$, namely  a horizontal shear layer with a small perturbation in the vertical velocity.
 The flow quickly degenerates into a periodic array of vortices, whose smallest scale is related to the laminar viscosity.
For this test,  one single periodic patch is used. Then, the splines of our (conformal) FEEC solver are set up with polynomial degree $\sf{p}=2$ and $100\times 100$ cells per patch. 
The contour plots for the vorticity $\omega=\curl \uu$ are presented  in figures \ref{fig:DSL1}-\ref{fig:DSL2} for the times $t=1.5$, $2.0$, $2.5$, $3.0$, $3.5$ and $4.0$ next to the reference solution computed  with a  5th order accurate spectral discontinuous Galerkin scheme on a staggered $40\times40$ grid, see \cite{FambriDumbser}.
Despite the presence of small oscillations, the overall dynamics is not corrupted. On the contrary, the numerical solution seems to match very well the proposed reference solution.

\section{Conclusion and perspectives}
\label{sec:Conclusion}
In this article we have proposed a new formulation of the Navier-Stokes equations using operators from the exterior calculus framework. By discretizing it with conforming and 
broken FEEC spaces we have next derived two numerical schemes which preserve the mass, momentum and energy of the fluid, in addition to being pressure robust. 
A middle point time discretization was then proposed to preserve these invariants at the fully discrete level and was analysed which resulted in a CFL condition for our scheme. By applying our numerical method to several standard test cases we have then verified its conservative properties, its high order accuracy and its ability to handle general boundary conditions. In particular, our tests demonstrate the ability of our scheme to correctly simulate flows dominated by pressure gradients, viscosity or advection terms.

An interesting direction for future work will be to extend this scheme to more complex models such as the compressible Euler equation. Further extension to non-linear MHD systems is a natural step, given the natural ability of FEEC operators to preserve the structure of Maxwell equations at the discrete level.  
\vfill
\section*{Acknowledgements}
The authors thank the developers and maintainers from the Psydac team, in particular Yaman G\"u\c{c}l\"u and Said Hadjout, for their efficient help in using the library for the purpose of this study. Inspiring discussions with Eric Sonnendr\"ucker are also warmly acknowledged.

\clearpage 

%
%
%
%

\bibliographystyle{plain} 
\bibliography{ins}

\end{document}